\PassOptionsToPackage{unicode}{hyperref}
\PassOptionsToPackage{hyphens}{url}
\PassOptionsToPackage{dvipsnames,svgnames,x11names}{xcolor}
\documentclass[12pt]{article}
\makeatletter
\ifx\paragraph\undefined\else
  \let\oldparagraph\paragraph
  \renewcommand{\paragraph}{
    \@ifstar
      \xxxParagraphStar
      \xxxParagraphNoStar
  }
  \newcommand{\xxxParagraphStar}[1]{\oldparagraph*{#1}\mbox{}}
  \newcommand{\xxxParagraphNoStar}[1]{\oldparagraph{#1}\mbox{}}
\fi
\ifx\subparagraph\undefined\else
  \let\oldsubparagraph\subparagraph
  \renewcommand{\subparagraph}{
    \@ifstar
      \xxxSubParagraphStar
      \xxxSubParagraphNoStar
  }
  \newcommand{\xxxSubParagraphStar}[1]{\oldsubparagraph*{#1}\mbox{}}
  \newcommand{\xxxSubParagraphNoStar}[1]{\oldsubparagraph{#1}\mbox{}}
\fi
\makeatother

\usepackage{longtable,booktabs,array}
\usepackage{calc} 
\usepackage{etoolbox}
\makeatletter
\patchcmd\longtable{\par}{\if@noskipsec\mbox{}\fi\par}{}{}
\makeatother
\IfFileExists{footnotehyper.sty}{\usepackage{footnotehyper}}{\usepackage{footnote}}
\makesavenoteenv{longtable}
\usepackage{graphicx}
\makeatletter
\def\maxwidth{\ifdim\Gin@nat@width>\linewidth\linewidth\else\Gin@nat@width\fi}
\def\maxheight{\ifdim\Gin@nat@height>\textheight\textheight\else\Gin@nat@height\fi}
\makeatother
\setkeys{Gin}{width=\maxwidth,height=\maxheight,keepaspectratio}
\makeatletter
\def\fps@figure{htbp}
\makeatother

\makeatletter
\@ifpackageloaded{caption}{}{\usepackage{caption}}
\AtBeginDocument{%
\ifdefined\contentsname
  \renewcommand*\contentsname{Table of contents}
\else
  \newcommand\contentsname{Table of contents}
\fi
\ifdefined\listfigurename
  \renewcommand*\listfigurename{List of Figures}
\else
  \newcommand\listfigurename{List of Figures}
\fi
\ifdefined\listtablename
  \renewcommand*\listtablename{List of Tables}
\else
  \newcommand\listtablename{List of Tables}
\fi
\ifdefined\figurename
  \renewcommand*\figurename{Figure}
\else
  \newcommand\figurename{Figure}
\fi
\ifdefined\tablename
  \renewcommand*\tablename{Table}
\else
  \newcommand\tablename{Table}
\fi
}
\@ifpackageloaded{float}{}{\usepackage{float}}
\floatstyle{ruled}
\@ifundefined{c@chapter}{\newfloat{codelisting}{h}{lop}}{\newfloat{codelisting}{h}{lop}[chapter]}
\floatname{codelisting}{Listing}

\makeatother
\makeatletter
\@ifpackageloaded{caption}{}{\usepackage{caption}}
\@ifpackageloaded{subcaption}{}{\usepackage{subcaption}}
\makeatother

\usepackage{iftex}
\ifLuaTeX
  \usepackage{selnolig}  
\fi
\usepackage{amsmath,amssymb}
\usepackage[]{natbib}
\usepackage{bookmark}

\IfFileExists{xurl.sty}{\usepackage{xurl}}{} 
\hypersetup{
  pdftitle={Title},
  pdfauthor={Author 1; Author 2},
  pdfkeywords={3 to 6 keywords, that do not appear in the title},
  colorlinks=true,
  linkcolor={blue},
  filecolor={Maroon},
  citecolor={Blue},
  urlcolor={Blue},
  pdfcreator={LaTeX via pandoc}}

\usepackage{tikz}
\usepackage{tikz-cd}
\usepackage{amsthm}
\usepackage{algorithm}
\usepackage{algpseudocode}

\graphicspath{{}}

\newtheoremstyle{thmstyleone}%
  {\topsep}{\topsep}%
  {\itshape}%
  {}%
  {\bfseries}%
  {}%
	{ }%
  {\thmname{#1}\thmnumber{ #2}\thmnote{ (#3)}}

\newtheoremstyle{thmstyletwo}%
  {\topsep}{\topsep}%
  {}%
  {}%
  {\bfseries}%
  {}%
	{ }%
  {\thmname{#1}\thmnumber{ #2}\thmnote{ (#3)}}

\newtheoremstyle{thmstylethree}%
  {\topsep}{\topsep}%
  {}%
  {}%
  {\bfseries}%
  {}%
	{ }%
  {\thmname{#1}\thmnumber{ #2}\thmnote{ (#3)}}

\theoremstyle{thmstyleone}%
\newtheorem{theo}{Theorem}
\newtheorem{propo}[theo]{Proposition}%
\newtheorem{lemma}[theo]{Lemma}%
\theoremstyle{thmstyletwo}%
\newtheorem{example}{Example}%
\newtheorem{rem}{Remark}%
\theoremstyle{thmstylethree}%
\newtheorem{defi}{Definition}
\newtheorem{assumption}{Assumption}
\usepackage[doublespacing]{setspace}
\usepackage[fontsize=12pt]{fontsize}

\usepackage{url}

\usepackage{mathtools}

\newcommand{\anon}{1}

\begin{document}

\def\spacingset#1{\renewcommand{\baselinestretch}%
{#1}\small\normalsize} \spacingset{1}


\if1\anon
{
  \title{\bf Causal inference via implied interventions}
  \author{Carlos G. Meixide$^{1,2}$\thanks{Work done at the University of California, Berkeley} \\ 
    Mark J. van der Laan$^{3}$ \\
    \phantom{s} \\
	$^{1}$Instituto de Ciencias Matemáticas, CSIC, Spain. \\
	$^{2}$Departamento de Matemáticas, Universidad Autónoma de Madrid, Spain. \\
	$^{3}$Division of Biostatistics, University of California, Berkeley, CA, USA.}
  \maketitle
} \fi

\if0\anon
{
  \bigskip
  \bigskip
  \bigskip
  \begin{center}
    {\LARGE\bf Causal Inference via Implied Interventions}
\end{center}
  \medskip
} \fi

\bigskip
\begin{abstract}
	In the context of having an instrumental variable, the standard practice in causal inference begins by targeting an effect of interest and proceeds by formulating assumptions enabling its identification. We turn this around by adhering to the interventions the observational distribution allows to identify, rather than starting with a desired causal estimand and imposing untestable conditions. The randomization of an instrument and its exclusion restriction define a class of auxiliary stochastic interventions on the treatment that are implied by stochastic interventions on the instrument. This mapping characterizes the identifiable causal effects of the treatment on the outcome given the observable distribution.  The identified effect is the impact of a stochastic encouragement by the instrument that propagates through the unaltered treatment selection mechanism, rather than the effect of a hypothetical intervention that overrides how treatment is naturally chosen. Alternatively, searching for an intervention on the instrument whose implied one best approximates a desired target naturally leads to a projection representing the closest identifiable treatment effect. The generality of this projection allows to select different norms and indexing functional sets that give rise to diverse estimation problems, some of which we address using Expectation-Maximization and the Highly Adaptive Lasso.
\end{abstract}

\noindent%
\vfill

\newpage
\spacingset{1.8} 

\section{Introduction}

The main problem of estimand identification under hidden confounding is the inability to use \citet{robins1986new}'s G-computation formula.
We assume throughout the paper the following nonparametric structural equation model \citep[NPSEM,][]{pearl2000causality}. \begin{assumption}\label{ass:natsem} $W = f_W(U_W), Z = f_Z(W, U_Z), A = f_A(Z, W, U_A), Y = f_Y(A,W,U_Y) $\end{assumption}

Each random variable $V$ is expressed as an arbitrary function of its parents $f_V$ together with an error term $U_V$. Given a distribution over the error terms, the equations then specify a joint distribution over the observable variables that we term the \textit{natural} distribution $P$. In this setup, $W$ denotes a vector of observed baseline covariates, $Z$ an instrumental variable satisfying

\begin{assumption}\label{ipr}(Conditional randomization): $  Z \perp U_A, U_Y |W,$\end{assumption}

\noindent $A$ is the treatment or exposure of interest and $Y$ the response or outcome variable. The primary goal is to evaluate causal effects of $A$ on $Y$, while accounting for potential unmeasured confounding.  {Note that we do not allow $Z$ to appear in the specification of $f_Y$, assumption commonly known as \textit{exclusion restriction}.} In contrast with \citet{pearl2000causality}'s original NPSEM formulation, we introduce the presence of hidden confounding by not specifying any conditional independence restrictions for the residuals in Assumption \ref{ass:natsem} except for the particular one in Assumption \ref{ipr} commonly referred to in the literature as \textit{instrument independent of hidden confounders}. A historical illustration is the return-to-education design of \citet{card1995using}: proximity to a college serves as the instrument $Z$, years of schooling as the treatment $A$, and earnings as the outcome $Y$, with proximity plausibly unrelated to the unobserved determinants of ability and wages encoded in $U_A$ and $U_Y$.

We write $P$ for the natural law. For any variables, $p(\cdot \mid \cdot)$ denotes conditional densities under $P$, with $g(A \mid W) := p(A \mid W)$ the treatment propensity score and $h(Z \mid W)$ the instrument density. Let $H(\cdot \mid W)$ be the conditional measure with density $h(\cdot \mid W)$ with respect to a dominating measure $\mu(\cdot \mid W)$; we use uppercase letters for measures and lowercase for densities, and write $H^*(\cdot \mid W)$ for the measure with density $h^*(\cdot \mid W)$. For an instrument policy $h^*$, define the induced treatment marginal $g(h^*)(A \mid W) := p(A^{h^*} \mid W)$. Theorem~\ref{idg} shows that
\[ g(h^*)(A \mid W) = \int p(A \mid Z, W)\, dH^*(Z \mid W). \]

To fix ideas, in the binary case $g(h^*)(w) = g_0(w)(1-h^*(w)) + g_1(w)h^*(w)$, where $g_z(w):=P(A=1\mid Z=z, W=w)$. A toy NPSEM and simulation that illustrate the map $h^* \mapsto g(h^*)$ are given in the Supplement.

{Figure~\ref{fig:gmap} (left) depicts the space of all possible stochastic interventions $h^*$, while Figure~\ref{fig:gmap} (right) shows their \textit{implied} interventions. Each intervention $h^*$ is assigned the same color as its image $g(h^*)$. The dot on the left represents the observational $h$, which is mapped to $g(h)$, the dot on the right, corresponding to the observational propensity score $P(A=1 \mid W=w)$.}

\begin{figure}[t!]
\centering
\includegraphics[width=0.65\textwidth]{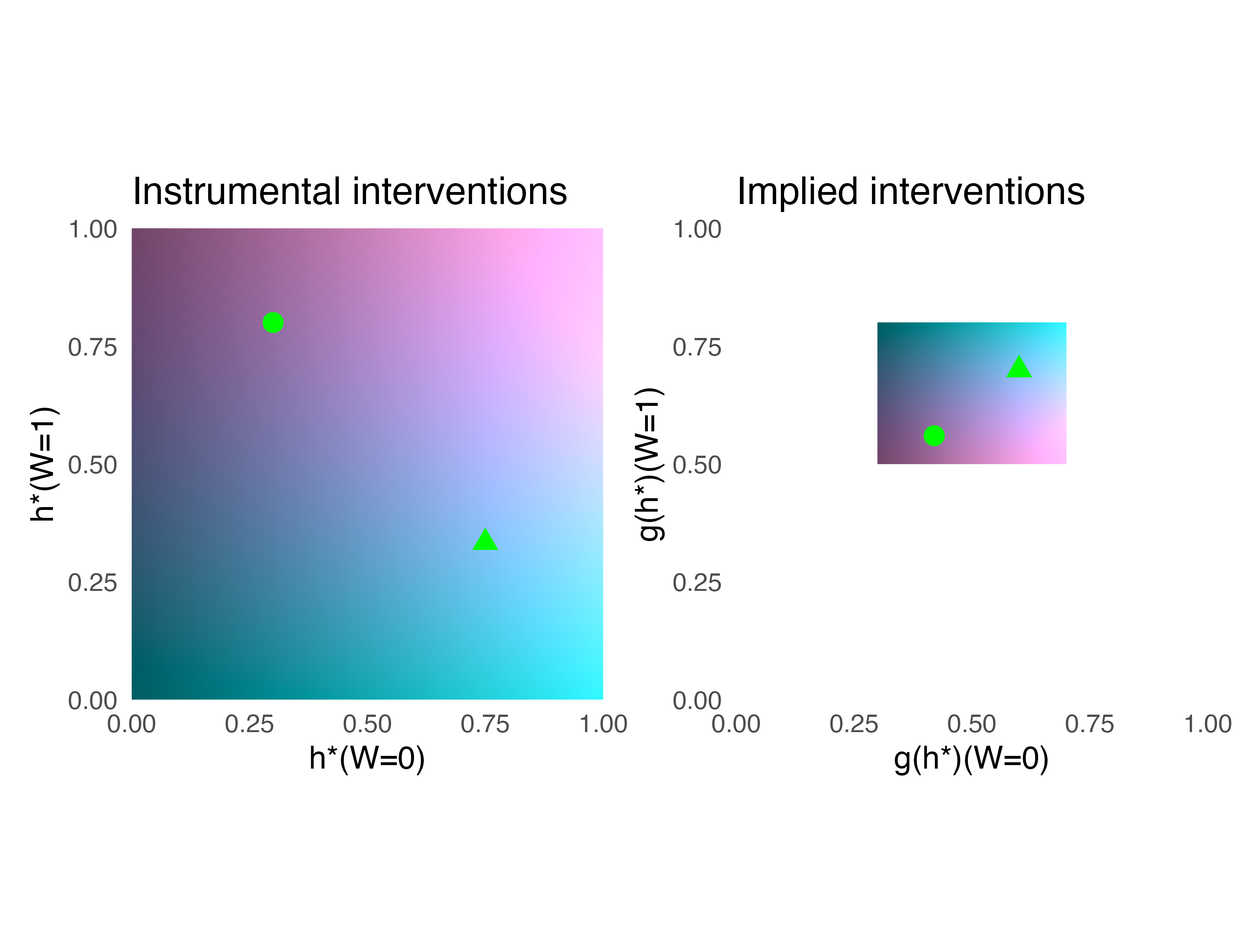}
\caption{Left: space of possible stochastic interventions $h^*$. Right: their images $g(h^*)$, with matching colors. Identification holds whenever $g^*$ has a preimage.}
\label{fig:gmap}
\end{figure}

\subsection{Literature review}

We now discuss two influential streams of the instrumental variables literature that are both compatible with our Assumption \ref{ass:natsem} and involve the observed variables $(W, Z, A, Y)$. While differing in historical development and primary applications, both frameworks rely on structural assumptions that map naturally to our causal model.

The first strand arises in the context of randomized experiments with noncompliance. In such settings, treatment assignment is randomized, but participants may not adhere to the assigned treatment. This leads to challenges in estimating the causal effect of the treatment actually received. A prominent solution, introduced by \citep{imbens} and further formalized by \citep{angrist1996identification}, is to combine Assumption \ref{ipr} together with an exclusion restriction and a monotonicity assumption to identify the local average treatment effect (LATE)---that is, the average causal effect for the subpopulation of compliers, individuals whose treatment status aligns with their assignment (see, e.g., \citep{tchetgen2024} for recent discussion on technical details). Importantly, this approach represents a shift in focus from the average treatment effect (ATE) in the entire population to an estimand defined over a latent subgroup. Since compliance status is never fully observed and therefore constituting a latent variable, LATE is not defined for identifiable units. Moreover, the LATE is inherently instrument-specific---its interpretation and value depend on the nature of the instrument and the group of compliers it induces. Different instruments may induce different compliance behaviors, and thus define different complier subgroups and corresponding LATEs. This limits the external validity of LATE: while it is a well-defined causal estimand, its scope is local, and generalizations beyond the induced complier group require further assumptions or alternative designs.

The second, more classical strand emerges from the econometric tradition of simultaneous equations, where the instrument helps resolve endogeneity arising from equilibrium relationships among economic variables.
Although non additive versions can also be found in the literature \citep{matzkin2003nonparametric, imbens2009identification}, this line of work often assumes $f_Y(A,W,U_Y) = f_P(A,W) + U_Y$, where $f_P$ solves the conditional moment equation $E [Y | Z]= E [f_P(A,W) | Z]$ derived from exclusion restriction $E [\varepsilon_Y | Z] = 0$. Solutions to this integral equation are typically ill-posed, owing to limited instrument strength relative to the treatment. While solutions may exist theoretically, they are numerically unstable \citep{darolles}. \citet{neweak} take a fundamentally different perspective on inference under weak identification by focusing on inference for linear functionals of the ill-posed nuisance, such as the NPIV regression function. Analogously, \citet{vanderlaan2025nonparametricinstrumentalvariableinference} introduce projection parameters for linear inverse problems that remain valid without identification, establishing asymptotic results for linear functionals with many weak instruments and proposing cross-fold debiased estimators. Please refer to the Supplement for a more detailed review of the problem. 

We take a different perspective to address the limitations imposed by unmeasured confounding and instrument weakness. Rather than formulating unverifiable assumptions, we acknowledge that the observational distribution $P$ imposes intrinsic constraints on which causal effects are identifiable. For those that are not, we focus on their closest identifiable ones in the sense of certain norms, framing the problem as an optimization task. 

Specifically, we propose to intervene stochastically on the instrument $Z$, which is assumed to be randomized conditional on covariates $W$. This induces a corresponding change in the distribution of $A$, effectively defining an implied stochastic intervention on the treatment. The key insight is that the effect of this induced intervention on $Y$ is equivalent to the effect of manipulating $Z$ in the specific manner that produces the observed perturbation in $A$. We then leverage this interventional auxiliary distribution of $A$ to compute the post-intervention expectation of $Y$ via a G-computation formula.

\section{Instrumental stochastic interventions}
The hardness of identifying causal effects of deterministic interventions opens the door to the more stable alternative of manipulating the assignment process itself instead of hard fixing values outright. For instance, identifying causal effects by deterministically fixing an exposure leans on assumptions rarely satisfied in real-world settings.  {Stochastic interventions were introduced in the unconfounded-treatment setting by \citet{munoz2012population,diaz2013assessing}, generalized to modified treatment policies that depend on the natural value of treatment by \citet{haneuse2013estimation,young2014identification}, extended in various directions by \citet{kenshift,diaz2023nonparametric,levis2024stochastic}, and applied to instrument-based designs by \citet{mauro,rakshit2024,Qiu02012021}.} These interventions are not only mathematically better-posed, but also more realistic from a clinical viewpoint, while still recovering classical deterministic interventions as a limiting case \citep{rakshit2024,mauro}.

In addition, it is typically the instrument---not the treatment---that is subject to intervention during experimentation. This reflects the reality of many applied settings where direct manipulation of the treatment $A$ is infeasible or unethical, while influencing the instrument $Z$ is both realistic and policy-relevant. For example, in encouragement designs, common in health and education studies, individuals are not forced to take a treatment but are instead encouraged through informational nudges or subsidies \citep{Zelen1979,Zelen1990}. Similarly, recommendation systems in online platforms modify the exposure $Z$ (e.g., what content or products a user sees) rather than the user's final decision $A$. 

Assume an intervened scenario in which the instrument $Z$ is generated not from its natural structural function $f_Z$, but instead from an alternative function $f_Z^*$.  {Crucially, we do \emph{not} introduce new versions of the endogenous error terms $U_A,U_Y$: the intervention acts on the $Z$-equation alone and the downstream equations for $A$ and $Y$ are propagated through the \emph{same} NPSEM noise as in the natural world. This keeps intact whatever dependence exists between $U_A$ and $U_Y$ under hidden confounding, and is precisely why our estimand is meaningful without destroying the confounding structure of the natural environment.}
\begin{assumption}\label{ass:intsem} $W = f_W(U_W) , Z^{h^*}= f^*_Z(W, {U}_Z) , A^{h^*} = f_A(Z^{h^*}, W,  {U_A}) , Y^{h^*} = f_Y(A^{h^*},W, {U_Y}) $\end{assumption}

Here $f^*_Z$ is the new structural equation encoding the new interventional conditional density $h^*$.  {We additionally require that the intervened instrument remains properly randomized in the new environment:}

\begin{assumption}[Intervened instrument properly randomized]\label{iipr} $Z^{h^*} \perp  {U_A, U_Y}  |W$
\end{assumption}

For instance, when the treatment is binary, \citep{kenshift} proposes in the unconfounded treatment setting a class of stochastic interventions  {that replace the natural propensity score $\pi(W) := P(Z=1\mid W)$ by a shifted version $\pi_\delta(W) = \frac{\delta\,\pi(W)}{\delta\,\pi(W) + 1 - \pi(W)}$, for $\delta>0$. This yields $\delta = \frac{\pi_\delta(W)/(1-\pi_\delta(W))}{\pi(W)/(1-\pi(W))}$, so that $\delta$ has the interpretation of an odds ratio.} They show that the effects of these interventions on the response are nonparametrically identified without requiring positivity.  {An extension to continuous treatments is developed by \citet{schindl2024incremental}. We adapt this construction to the instrument $Z$ in our framework, where the lever is the policy on $Z$ rather than on the treatment itself.}

So far, a new conditional density $h^*(\cdot \mid W)$ is introduced, replacing the original $h(\cdot \mid W)$, thereby modifying the distribution of $Z$ given covariates $W$. One may conceptualize the intervened NPSEM in Assumption \ref{ass:intsem} as a new environment arising from manipulating the conditional distribution of the instrument $Z \mid W$. This perspective is consistent with the practical view of instruments as variables that can, at least in principle, be manipulated by the analyst or policy maker. 

This setup reflects a stochastic intervention on $Z$, where rather than deterministically setting $Z$ to a fixed value (as in a ``hard'' intervention), we intervene by changing the probability distribution from which $Z$ is drawn. The idea of modifying the stochastic mechanism behind $Z$ is formalized in \citep{Qiu02012021} under the name of an individualized encouragement rule. The key insight is that such rules allow for a softer, more flexible form of intervention---one that remains within the observational data distribution's support and avoids the practical and theoretical limitations of deterministic assignments. 

Assumption~\ref{wp}  {constitutes} a technical guarantee allowing to use the tower property of conditional expectations and TMLE in Section~\ref{sec:semi}.

\begin{assumption}[Instrument positivity]\label{wp} $\frac{h^*}{h} <\infty$ a.e.
\end{assumption}

Note that Assumption \ref{wp} implies $h^* \ll h$, meaning $h^*$ absolutely continuous with respect to $h$: if $h(z|W)=0$ for some $z$ in the support, then $h^* (z|W)=0$.

The parameter $EY^{h^*}$ represents a causal quantity because it is in principle a functional of the post interventional distribution of $Y$ resulting from the replacement $f_Z \to f^*_Z$ in the original NPSEM in Assumption \ref{ass:natsem}, quantifying a causal effect of $Z$ on $Y$. Therefore, this expectation is not a purely associational quantity; it explicitly depends on the modified $Z|W$ mechanism. It answers a ``what if'' question: ``What would the expected value of $Y$ be if given $W$ we forced $Z$ to follow $h^*(Z|W)$?''. The fact that $EY^{h^*}$  itself is counterfactual does not have anything to do with whether it is identifiable or not. 

\begin{theo}[Identifiability of $EY^{h^*}$]\label{idth}
Under Assumptions \ref{ass:natsem}, \ref{ipr}, \ref{ass:intsem}, \ref{iipr} and \ref{wp}: $E[Y^{h^*}] = E_P E_{Z | W \sim h^*} \left[ E_P[Y \mid Z, W] |W \right]$. Equivalently $E[Y^{h^*}] = \int \left( \int  E_P[Y \mid Z, W] dH^*(Z|W)\right)dP(W)$
\end{theo}
 {Assumption~\ref{wp} ensures that the inner integral against $H^*(\cdot\mid W)$ is well-defined: whenever the inner conditional expectation $E_P[Y\mid Z,W]$ is undefined at a point $(z,w)$ because $h(z\mid w)=0$, absolute continuity $h^*\ll h$ forces $h^*(z\mid w)=0$ at the same point, so that point contributes zero to the integral. If $Z$ were also allowed to enter the structural function $f_Y$, Theorem~\ref{idth} would still hold; we exclude this case because the identifiability of the implied intervention in the next section relies on exclusion restriction.}

Theorem \ref{idth} is an actual identifiability result because the conditional expectations therein involve known mechanisms, which are from inner to outer: the natural conditional pdf of the response given the instrument and the covariates, the intervened $h^*$ (chosen by the user) and the marginal pdf of the covariates.

\begin{figure}[t!]
\centering
\begin{tikzcd}[column sep=large, row sep=large]
	& Z \arrow[rr, "\textcolor{blue}{EY}"] \arrow[dr, swap, "\textcolor{blue}{g(h)}"] & & Y \\
	& & A\arrow[ur, dashed, ->, "\quad  {EY}"{pos=0.25, anchor=west}]
\end{tikzcd}
\hspace{1.5em}
\begin{tikzcd}[column sep=large, row sep=large]
	& Z^{h^*} \arrow[rr, "\textcolor{blue}{EY^{h^*}}"] \arrow[dr, swap, "\textcolor{blue}{g(h^*)}"] & & Y^{h^*} \\
	& & A^{h^*} \arrow[ur, dashed, ->, "\quad  {EY^{g(h^*)}}"{pos=0.25, anchor=west}]
\end{tikzcd}
\caption{The blue estimands are identified due to the proper randomization of the instrument.  {Exclusion restriction implies that $Y^{h^*}$ can be viewed as mediated through $A^{h^*}$, whose induced conditional marginal given $W$ is $g(h^*)$. The dashed arrow labelled $EY^{g(h^*)}$ denotes the same counterfactual mean in the NPSEM-propagated world, \emph{not} the effect of drawing $A$ independently from $g(h^*)$; see Remark~S1 in the Supplementary Material.}}\label{fig:envs}

\end{figure}

\section{Implied stochastic interventions}\label{sec:implied}

The structural equation for $A^{h^*}$ in Assumption \ref{ass:intsem} reveals that  {modifying} the conditional distribution of $Z$ given $W$  {also modifies} the conditional distribution of $A$ given $W$.  {The marginal distribution of $A^{h^*}\mid W$ that emerges in the intervened world is thus a by-product of the policy on $Z$, not a free-standing policy density that one decides to impose on the treatment. This distinction is central to interpreting our results: throughout the paper the symbol $g(h^*)(A\mid W)$ is a \emph{descriptive summary} of what the intervention on $Z$ does to the distribution of $A$ inside the intervened world, and \emph{not} the density of a separate stochastic policy on $A$ drawn independently of the endogenous outcome noise $U_Y$.}

 {To fix ideas, consider the alternative, fully artificial world in which the treatment is assigned by drawing $\widetilde A^{g^*} \sim g^*(\cdot\mid W)$ \emph{independently} of $U_Y$ and then $\widetilde Y^{g^*} = f_Y(\widetilde A^{g^*}, W, U_Y)$. Even when $g^* = g(h^*)$, in general $E\!\left[Y^{h^*}\right] \;\neq\; E\!\left[\widetilde Y^{g^*}\right]$, because in the first world $A^{h^*} = f_A(Z^{h^*},W,U_A)$ still depends on $U_A$, which under hidden confounding is correlated with $U_Y$, whereas the independent draw breaks that dependence. Precise conditions under which the two counterfactuals agree are given in Remark~S1 of the Supplementary Material.}

\begin{theo}\label{propoim}
 {Let $h^*$ be a stochastic intervention on the instrument. Under Assumptions \ref{ass:natsem}, \ref{ipr}, \ref{ass:intsem} and \ref{iipr}, the counterfactual outcome under the intervention on $Z$ satisfies $Y^{h^*} \;=\; f_Y\!\left(A^{h^*},\,W,\,U_Y\right)$,
so that $Y^{h^*}$ depends on $Z^{h^*}$ only through $A^{h^*}$. Consequently, $E[Y^{h^*}]$ is a functional of $P$ and $h^*$ whose induced conditional marginal on $A$ given $W$ is $g(h^*)(A\mid W) := p(A^{h^*}\mid W)$. In particular, writing $E[Y^{g(h^*)}]$ for the expectation in this NPSEM-propagated world, one has $E[Y^{h^*}] = E[Y^{g(h^*)}]$.}
	\end{theo}
	
	 {\textit{Proof.} Since $Z$ does not appear in $f_Y$ (exclusion restriction, built into Assumption \ref{ass:natsem}), $Y^{h^*} = f_Y(f_A(Z^{h^*},W,U_A),W,U_Y) = f_Y(A^{h^*},W,U_Y)$. By construction $A^{h^*}\mid W \sim g(h^*)(\cdot\mid W)$, hence any functional of the joint law of $(Y^{h^*},W)$ --- in particular its mean --- is a functional of the triple $(g(h^*), p(Y\mid A,Z,W), P_W)$ alone, which justifies the notation $E[Y^{g(h^*)}]$ for $E[Y^{h^*}]$. Note that no cross-world equality between different versions of $U_Y$ is required: the original NPSEM noise is kept throughout.}  {Positivity (Assumption~\ref{wp}) is not invoked: the statement is a structural identity between counterfactual random variables and does not involve an integral against $h^*/h$.}
	
	 {The reader should be careful: the notation $E[Y^{g(h^*)}]$ is a relabelling of $E[Y^{h^*}]$ in the NPSEM-propagated world, \emph{not} the expectation of the outcome under an independent draw of $A$ from $g(h^*)$. The two coincide only under the conditions of Remark~S1 in the Supplementary Material.}
	
	\begin{rem}[Equivalent derivation without an NPSEM]\label{rem:po}
 {The identification formula in Theorem~\ref{idth} can also be derived directly in the potential-outcomes language, without committing to an NPSEM. Let $\{Y^a : a\in\mathcal A\}$ and $\{A^z : z\in\mathcal Z\}$ be potential outcomes for $Y$ and $A$, and assume the IV core together with exclusion restriction, phrased as the conditional independence statement $(Y^a, A^z) \perp Z \mid W$ for all $(a,z)\in\mathcal A\times\mathcal Z$.}  {Here $Y^a$ is indexed by the treatment value, not the instrument: exclusion is what makes $Y^a$ well defined, since $Y^{z,a}=Y^a$. The G-computation formula below recovers the cross-world joint $p(Y^a,A^z\mid W)$ from $p(Y,A\mid Z=z,W)$, and this step requires the joint statement $(Y^a,A^z)\perp Z\mid W$; the reduced-form independence $Y^z\perp Z\mid W$, although implied, does not by itself identify that joint.}

 {Then for any stochastic intervention $h^*$ on $Z$,
	\begin{equation}\label{eq:po-id}
		E\!\left[Y^{h^*}\right] \;=\; \int\!\!\int\!\!\int y\; p\!\left(Y^a=y,\, A^z=a \,\middle|\, W\right)\, dH^*(z\mid W)\, d\mu(a)\, d\mu_Y(y)\, dP(W).
	\end{equation}
	By contrast, a physical intervention that draws $\widetilde A \sim g(h^*)(\cdot\mid W)$ \emph{independently} of $Y^a$ yields
	\begin{equation}\label{eq:po-clean}
		E\!\left[\widetilde Y^{g(h^*)}\right] \;=\; \int\!\!\int\!\!\int y\; p\!\left(Y^a=y \,\middle|\, W\right)\, g(h^*)(a\mid W)\, d\mu(a)\, d\mu_Y(y)\, dP(W).
	\end{equation}
	Because of hidden confounding, $p(Y^a, A^z\mid W) \neq p(Y^a\mid W)\, p(A^z\mid W)$ in general, so \eqref{eq:po-id} and \eqref{eq:po-clean} differ. They coincide under (i) or (ii) of Remark~S1 in the Supplementary Material: under (i), the factorization $p(Y^a,A^z\mid W) = p(Y^a\mid W) p(A^z\mid W)$ holds directly; under (ii), integrating $y$ against $m(a,W) + v(W,U_Y)$ returns a quantity depending only on the marginal of $A^z$ given $W$. Our identification formula targets~\eqref{eq:po-id}: the intention-to-treat-style effect that can be recovered from observed data under instrumental randomization alone.}
\end{rem}
	
	 { The choice between \eqref{eq:po-id} and \eqref{eq:po-clean} is dictated by how treatment is generated in the system under study, not by analytic convenience. In clinical practice, the treatment $A$ a patient receives is the output of a decision process informed by latent prognostic features that are not in $W$: in intensive care, for example, treatment choices are shaped by acuity, frailty, and anticipated outcome, none of which are recorded. The same unobservables drive $Y$, so $U_A$ and $U_Y$ are conditionally dependent given $W$. This dependence is not a defect of the data; it is the mechanism by which treatment is actually allocated. The ICH E9(R1) framework codifies the corresponding analytic stance via the \emph{treatment-policy strategy}: post-randomization events such as treatment discontinuation or use of rescue medication are not hypothetically eliminated but allowed to propagate into the estimand, because they are part of how the intervention behaves once deployed \citep{ich2019e9r1,wei2025estimands}.}

 {Our framework adopts the same stance at the structural level. Under Assumption \ref{ass:intsem}, a stochastic policy on the instrument induces $A^{h^*}$ and $Y^{h^*}$ through the unaltered $f_A$ and $f_Y$ with the natural errors $U_A, U_Y$ retained, and identifies \eqref{eq:po-id}. The joint $p(Y^a, A^z \mid W)$, which carries the confounding between treatment and outcome, is preserved by construction. The atomic policy \eqref{eq:po-clean}, by contrast, factorizes this joint and corresponds to a thought experiment in which the analyst overrides the selection mechanism and draws $A$ from $g^*$ independently of $U_Y$. Such an atomic policy is rarely realizable when treatment decisions are informed by features shared with the outcome, and its identification would require either no hidden confounding or restrictions on $f_Y$ that the IV core does not impose. The induced marginal $g(h^*)(A \mid W)$ is therefore a \emph{description} of how the treatment distributes itself inside the encouragement world, not a free-standing policy on $A$. The question we answer is what happens to $Y$ when the antecedent encouragement mechanism is shifted by $h^*$, while the rest of the system---including the dependence between $A$ and the unmeasured drivers of $Y$---is left intact.}

	 {Theorem \ref{propoim} and the two preceding remarks together say that} the effect of an instrumental stochastic intervention identified by Theorem \ref{idth} is actually a treatment effect,  {mediated solely by changes in $A^{h^*}$, whose conditional marginal given $W$ is captured by $g(h^*)$}.
	\begin{defi}
Let $h^*$ be a stochastic intervention on the instrument.  {We call $p(A^{h^*}\mid W)=: g(h^*)(A\mid W)$ the \emph{induced} (or \emph{implied}) conditional marginal of $A$ in the world where $Z$ is drawn from $h^*(\cdot\mid W)$ and $A$ is generated through $f_A$ with the natural $U_A$. It is a description of $A$'s distribution inside that world, not a free-standing policy on $A$ to be interpreted in isolation.}
\end{defi}

Up to now, we have not made any assertion on whether it is possible to identify $g(h^*)$.  {The same rationale used in Theorem~\ref{idth} yields the following identification result.}

\begin{theo}[Identifiability of the implied intervention]\label{idg}
Under Assumptions  \ref{ass:natsem}, \ref{ipr},  \ref{ass:intsem} and \ref{iipr}: $g(h^*)(A|W) =
E_{Z \mid W \sim h^*}
p(A \mid Z, W)$
or equivalently
$g(h^*)(A|W)  = \int   p(A \mid Z, W) dH^*(Z|W)$

\end{theo}

 {\textit{Proof.} By construction $A^{h^*} = f_A(Z^{h^*}, W, U_A)$, so the conditional law of $A^{h^*}$ given $W$ is
$P(A^{h^*}=a \mid W) = \int P(f_A(z,W,U_A) = a \mid Z^{h^*}=z, W)\, dH^*(z\mid W)$. Under Assumption~\ref{iipr}, $Z^{h^*} \perp U_A \mid W$, so $P(f_A(z,W,U_A)=a \mid Z^{h^*}=z, W) = P(f_A(z,W,U_A)=a \mid W) = p(A=a\mid Z=z, W)$, where the last equality uses Assumption~\ref{ipr} (the natural $Z \perp U_A \mid W$) to identify the structural conditional from the observed data. Substituting and noting that positivity (Assumption~\ref{wp}) is not required here---the integral is over a probability measure $H^*(\cdot\mid W)$ and $p(A\mid Z,W)$ is bounded---yields the stated identity. \hfill$\square$}

This means that identifiability of the causal effect on the response after manipulating the instrument (guaranteed by Theorem \ref{idth}) comes for free with identification of  {the induced marginal $g(h^*)$ of $A$}.  {Because Theorem~\ref{idg} identifies $g(h^*)$ from the observed data, and Theorem~\ref{propoim} shows that the effect of the instrument on $Y$ is mediated entirely through that induced marginal, the causal effect of any stochastic policy on $Z$ is accessible without any further assumptions.}

We have arrived to a mapping that is of crucial importance for our work \begin{equation}\label{main:operator} 
H^*(Z|W) \overset{g}{\longmapsto} \int   p(A \mid Z, W) dH^*(Z|W) \end{equation}

 {The range of this operator characterizes all treatment marginal distributions that can be induced by policies on $Z$, and whose associated causal effects on $Y$ are therefore identifiable. The only aspect of the observed distribution that enters is the conditional assignment $p(A\mid Z,W)$.}

 {In general, for an arbitrary target treatment distribution $g^*$, the counterfactual mean $E[Y^{g^*}]$ under an atomic policy that draws $A$ from $g^*$ independently of $U_Y$ cannot be identified from the observed data unless $A$ is unconfounded given $W$ (see Remark~S1 in the Supplementary Material). By contrast, for any stochastic intervention $h^*$ on the properly randomized instrument $Z$, the counterfactual mean $E[Y^{h^*}]$ remains identifiable through Theorem~\ref{idth}, and Theorem~\ref{idg} reveals that $h^*$ induces a treatment marginal $g(h^*)$ describing how $A$ distributes itself in that world. Theorem~\ref{propoim} guarantees $E[Y^{h^*}] = E[Y^{g(h^*)}]$, where the right-hand side is the causal effect propagated through the unaltered selection mechanism, not an atomic draw from $g(h^*)$.}

A direct application of Bayes theorem leads to the last result in this part of the article, which will be relevant for the development of Section \ref{sec:ls}. 

\begin{propo}\label{bayes}Under Assumptions \ref{ipr}, \ref{iipr} and \ref{wp}: $g(h^*)(A|W)= g(A|W)E\left[ \frac{h^*(Z|W)}{h(Z|W)} \mid A,W\right]$
\end{propo}

Proposition \ref{bayes} means that under assumptions \ref{ipr}, \ref{iipr} and \ref{wp} every implied intervention is of the form in its rhs. This is interesting since for each $W$ we multiply $g(A|W)$ by a factor (either shrinking or increasing) that equals the  {conditional} mean of $h^*/h(Z|W)$ given $A$ (all conditional on $W$).  Generally speaking, multiplicative stochastic interventions are of interest as they preserve power while still answering interesting questions.

 {Our framework requires no preliminary assumption on the strength of the instrument. When the instrument is weak in the sense that $p(A\mid Z,W)\approx p(A\mid W)$, the integral $\int p(A\mid Z,W)\,dH^*(Z\mid W) \approx p(A\mid W)$, so that $g(h^*)$ remains close to the observational propensity score $p(A\mid W)$ and the set of attainable induced marginals collapses around this value (as Proposition~\ref{bayes} also shows). The range of $g(\cdot)$ therefore expands with the strength of the instrument; in the degenerate case in which $A$ is independent of $Z$ given $W$, the range reduces to the singleton $\{p(A\mid W)\}$.}

\section{Randomized experiments under noncompliance}\label{sec:late}

Recent developments in the causal inference literature have advocated moving beyond LATE toward estimating marginal causal effects that hold broader policy relevance \citep{carneiro,mogstad2018using}. These  {effects} describe the average impact of hypothetical interventions across an entire target population rather than conditioning exclusively on latent subpopulations such as compliers. Consequently, estimating marginal effects better addresses policymakers' need for generalizable insights when evaluating large-scale expansions. 

For binary $A$ and $Z$, Theorem \ref{idg} establishes $g(h^*)(w) = g_0(w)(1-h^*(w)) + g_1(w)h^*(w)$. Suppose now that we fix $g^*$ and require  $g_0(w) \leq g^*(w) \leq g_1(w) $ when $g_1(w) > g_0(w) $ and 
$g_1(w) \leq g^*(w) \leq g_0(w) $ when  $g_1(w)<   g_0(w) $. Then defining 
\begin{equation}\label{hthatimplies}
h^*(w) := \frac{g^*(w) - g_0(w)}{g_1(w) - g_0(w)}
\end{equation}

we have that $0 \leq h^*(w) \leq 1$ and $g^*=g(h^*)$ ($h^*$ implies $g^*$). Observe that the requirements on $g^*$ resemble a generalized notion of monotonicity. The classical formulation of monotonicity (which we hereby refer to as \textit{hard} monotonicity) assumes $P(A^{0} \leq A^{1})  =1$; where $A^z = f_A(z,W,U_A)$ are potential outcomes obtained straight from fixing $Z=z$ in the structural equation for $A$ in the NPSEM in Assumption \ref{ass:natsem}. Note that under hard monotonicity $g_1(w) - g_0(w) = P(\textrm{complier}|W=w)$, which coincides with the denominator in the conditional Wald ratio: \citep{angrist1996identification, ogburn2015doubly} $$ \operatorname{Wald}(W):=\frac{E[Y| Z=1,W]- E[Y|Z=0,W]}{E[A| Z=1,W]- E[A|Z=0,W]}$$ 
We also point that the stronger the instrument, the wider range of implied interventions we can identify by inspection of the sandwich enclosing $g^*(w)$ present in the above conditions.

\begin{defi}[{$Z$-compatibility}]\label{def:weakmono}Let $Z$ and $A$ be binary and let $g_z(w) = P(A=1|Z=z,W=w)$.  {We say that a target treatment distribution $g^*(w)$ is \emph{$Z$-compatible} (also referred to as satisfying weak monotonicity in the earlier literature) when} $g^*(w) \in \left[ \min\{g_0(w), g_1(w)\}, \max\{g_0(w), g_1(w)\} \right]$, for $w$ in the support of $W$.  {We note that this condition constrains the desired policy $g^*$ rather than the potential outcomes $A^z$; it is not a weakening of the almost-sure monotonicity of \citet{imbens,angrist1996identification}, but rather a compatibility requirement between $g^*$ and the observed conditional distribution of $A$ given $Z,W$. When $g_0(w)=g_1(w)$ (weak instrument on the event $\{W=w\}$), $Z$-compatibility forces $g^*(w)=g_0(w)=g_1(w)$, so that the product $(g^*(w)-f^*(w))\operatorname{Wald}(W)=0$ in Theorem~\ref{generaldiff} below is defined as zero on this event.}

\end{defi}

 {Therefore, $Z$-compatibility is a condition on the target treatment distribution, not on the potential outcomes $A^z$, and it is empirically checkable given estimates of $g_0(w)$ and $g_1(w)$.} The next theorem shows that, even without enforcing $A^{(0)}\le A^{(1)}$ $P$-a.s., one can identify a broad class of population level causal contrasts $E[Y^{g^*}-Y^{f^*}]$ as long as both $g^*$ and $f^*$ are $Z$-compatible. An auxiliary lemma and the proof can be found in the Supplement.

\begin{theo}\label{generaldiff}
Let $Z$ and $A$ be binary and $f^*$, $g^*$ be any  {$Z$-compatible} target treatment distributions. Under the hypotheses of Theorem \ref{idth},
$$E\left[Y^{g^*} - Y^{f^*} \right]  = E\left[ \left(g^*(W) - f^*(W)\right)\cdot \operatorname{Wald}(W) \right]$$

If in addition $P$ satisfies hard monotonicity,  {denoting by $\operatorname{LATE}(W) := E[Y^1 - Y^0 \mid W, \mathrm{complier}]$ the conditional local average treatment effect on the subpopulation of compliers,}

$$E\left[Y^{g^*} - Y^{f^*} \right]  = E\left[ \left(g^*(W) - f^*(W)\right)\cdot \operatorname{LATE}(W) \right]$$

\end{theo}

 {On the event $\{W : g_0(W)=g_1(W)\}$, the $Z$-compatibility condition forces $g^*(W)=f^*(W)=g_0(W)$, so the product $\bigl(g^*(W)-f^*(W)\bigr)\operatorname{Wald}(W)$ is defined to be zero regardless of whether $\operatorname{Wald}(W)=g_1(W)-g_0(W)$ is zero or not. This convention avoids a $0/0$ indeterminacy and is consistent with the fact that a weak instrument on that event carries no identifying information; the integral over the complement of that event is well-defined without any relevance assumption.}

 {This result shows that identifiability of population-level causal contrasts is not restricted to \emph{local} effects on the complier subpopulation: under imperfect compliance, the best marginal contrast we can identify is the estimand on the left-hand side. Rather than requiring that hard monotonicity holds for every unit, it suffices to require that the target distributions $g^*$ and $f^*$ satisfy a strictly weaker, and empirically checkable, $Z$-compatibility condition.}

It is important to note that the hard monotonicity assumption typically invoked in the LATE literature is much stronger than simply requiring that \( Z = 1 \) increases the probability of receiving treatment \( A = 1 \). Rather, it is an almost surely statement: \(  {P}(A^{0} \leq A^{1}) = 1 \), which we would still need to identify LATE$:=E(Y^1-Y^0|\operatorname{compliers})$. 

Exact identification of $E(Y^1-Y^0)$ would require full compliance ($P(A=1\mid Z=0,W)=0$ and $P(A=1\mid Z=1,W)=1$)---a setting where IV methods are unnecessary. LATE is instead defined over the latent complier subpopulation ($A^1=1, A^0=0$). By contrast, our approach always identifies a marginal, population-wide contrast; the identified contrast narrows as the instrument weakens but remains well-defined and supported by the data.

\section{Semiparametric efficiency}\label{sec:semi}

Under a model that leaves the components of the joint distribution of $O = (W, Z, A, Y)$ unrestricted except for $h(Z \mid W)$, the efficient influence curve (EIC) for the parameter $\Psi_{h^*}(P)=E_P(Y^{h^*})$ characterizes the lowest possible asymptotic variance achievable by any regular and asymptotically linear estimator. According to efficiency theory \citep{bickel1993efficient} an estimator that solves the empirical EIC estimating equation for a known $h^*$ is asymptotically efficient in this model. Targeted Minimum Loss-based Estimation (TMLE) introduced by \citep{vandrlrubin2006} provides a principled approach to constructing such estimators. TMLE proceeds in two stages: first, it defines an initial estimator of the relevant parts of the data distribution. Then, it constructs a parametric submodel through this initial estimate whose score, at zero fluctuation, spans the EIC. The likelihood within this submodel is iteratively maximized until convergence, at which point the updated estimator satisfies the efficient influence curve estimating equation. The plug-in estimate of the target parameter based on this updated distribution is the TMLE. 

Our parameter of interest is $\Psi_{h^*}(P) = E_P(Y^{h^*})$, which we were able to identify in Theorem \ref{idth}, so that

\[
\begin{array}{rcl}
	\Psi_{h^*} : 	\mathcal{M} & \longrightarrow & \mathbb{R} \\
	P & \mapsto & E_P Y^{h^*}
\end{array}
\]

with $\mathcal{M}$ the fully unrestricted nonparametric model. 

\begin{propo}\label{prop:eic}
	The efficient influence curve of $\Psi_{h^*}(P) = E_P(Y^{h^*})$ is

	$$D^*_{h^*,P} = \frac{h^*(Z|W)}{h(Z|W)}(Y - E(Y|W,Z)) + \int E(Y|W,Z) h^*(Z|W)dZ - \Psi_{h^*}(P) $$

	 {This functional form is analogous to the EIC for stochastic interventions on an unconfounded treatment derived in \citet{munoz2012population} and \citet{diaz2013assessing}; see also \citet{diaz2020causal}. The result holds here because Theorem~\ref{idth} identifies $E[Y^{h^*}]$ through a G-computation formula in the observed-data model in which $Z$ plays the role of an unconfounded exposure, so the parameter $\Psi_{h^*}$ inherits the efficiency bound of that literature directly.}
\end{propo}

This is the efficient influence curve assuming a model on $h(Z|W)$ but leaves other factors of the joint density of $O=(W,Z,A,Y)$ nonparametric. We begin by constructing an initial estimator $Q_n^0(Z, W)$ of the outcome regression $E[Y \mid Z, W]$ and defining the \textit{clever weights} as:
$H(Z, W) := \frac{h^*(Z \mid W)}{h(Z \mid W)}$. In this case, the clever covariate is simply 1, and we incorporate $H(Z, W)$ as sampling weights in the targeting step. To update $Q_n^0$, we fit a fluctuation model using weighted logistic regression: \verb|glm(Y ~ offset(Qn0(Z, W)) + 1,  weights = h*(Z, W) / h(Z, W))|. This model treats $Q_n^0(Z, W)$ as a fixed offset and fits a single intercept parameter $\hat{\epsilon}$ using weights $H(Z, W)$. The estimated $\hat{\epsilon}$ adjusts $Q_n^0$ in the direction of the EIC. We update the initial estimator as $Q_n^{\text{updated}}(Z, W) = Q_n^0(Z, W)+ \hat{\epsilon}$. {For the target parameter $E[Y^{h^*}]$ with fixed $h^*$, we require only a single update step of the outcome regression. The procedure becomes iterative only when the target parameter is such that the clever weights themselves depend on the outcome regression, thereby requiring successive updates.}  In addition, we do not need to target the empirical distribution of $W$ since we always solve that component. 

\subsection{Observed or hidden confounders?}\label{hidnohid}

If the instrument is independent of the covariates---i.e., if the natural policy satisfies $h(Z \mid W) \;=\; h(Z)\,$,
then any subset of the observed covariates in $W$ can be treated as though it were unobserved for the purpose of identification of the causal effect on $Y$ of an intervened policy $h^*$ that involves covariates in the complement of such a subset. In particular, for each subset $S \subseteq \{1,\ldots,p\}$ corresponding to reduced data structures $W_S = (W_j : j \in S)$, one obtains the same G-computation formula:
\begin{equation}\label{eq:identification_family}
	E[Y^{h^*}] = \int \left( \int  E_P[Y \mid Z, W_S] dH^*(Z|W_{S'})\right)dP(W_S)
\end{equation}

for a chosen $H^*(Z|W_{S'})$ as far as the same assumptions in Theorem \ref{idth} are satisfied and $S' \subseteq S$. Observe that, although $W_S$ appears in the inner conditional expectation $E_P[Y \mid Z,\,W_S]$, the interventional policy $h^*(Z \mid W_{S'})$ need not depend on all components in $W_S$. Consequently, for any fixed $h^*(Z \mid W_{S'})$, the same causal estimand $E[Y^{h^*}]$ becomes identified using all the different subsets $S$ of the observed covariates containing $S'$ ---one expression for each $S \supseteq S'$. 

Furthermore, Theorem \ref{idg} produces  {a corresponding family of induced treatment marginals} for a fixed $h^*(Z \mid W_{S'})$:

\begin{equation}\label{eq:implied_interventions}
	g(h^*)(A|W_S)  = \int   p(A \mid Z, W_S) dH^*(Z|W_{S'}), \quad S \supseteq S'
\end{equation}

each of which yields exactly the same causal effect on the response.  {This multiplicity is useful when selecting an implied intervention for reporting: one can pick the induced marginal whose treatment distribution profile is simplest to describe or closest to a target $g^*$ of scientific interest, while noting that all members of the family identify the same effect.}

Finally, once we have agreed the $h^*(Z \mid W_{S'})$, one must use the full $W$ to build the EIC for TMLE in Proposition \ref{prop:eic} to maximize statistical efficiency. In this way, one makes the most efficient use of all available information while still maintaining parsimony in subsequent analyses and presentations by choosing a simple implied intervention to describe. We will adhere to this procedure in our own analysis of the 2008 Oregon health insurance experiment in Section \ref{sec:oregon}.

\section{Projecting treatment interventions}
Until now we discussed how a stochastically intervened instrument $Z^{h^*}$ leads to a counterfactual treatment $A^{h^*}$ which itself leads to an intervened response $Y^{h^*}=f_Y( A^{h^*},W, {U_Y})$ {, with $A^{h^*}$ and $U_Y$ generally dependent through the hidden confounder}.  {By contrast, drawing $\widetilde A^{g^*}\sim g^*(\cdot\mid W)$ independently of $U_Y$ leads to $\widetilde Y^{g^*}=f_Y(\widetilde A^{g^*},W,U_Y)$, and as discussed in Remark~S1 of the Supplementary Material this second world yields a different expectation in general.} Given an $h^*$, Theorem \ref{idg} shows how  {to} compute the $g^*=g(h^*)$  {describing the induced conditional marginal of $A^{h^*}$ given $W$}. 

 {This raises the following question: given a desired treatment-distribution profile $g^*$, does an instrument policy $h^*$ exist whose induced marginal matches $g^*$ exactly? If not, we seek the $h^\dagger$ whose induced marginal is closest to $g^*$ in a suitable sense. The set of all attainable induced marginals is $\mathcal{J} = \{ g(h^*) : h^* \in \mathcal{I}\}$, where $\mathcal{I}$ is the convex set of instrument-policy densities. Note that $g^*$ here is a \emph{reference} distribution summarizing how a researcher imagines treatment might be distributed under some hypothetical regime; the question is which instrument encouragement can bring the natural treatment-selection mechanism closest to that profile, while keeping the confounding structure intact. In the language of target trial emulation \citep{hernan2016target}, $\mathcal{J}$ characterizes the set of RCT-like treatment protocols that IV data can actually support, and $h^\dagger$ is the instrument policy that emulates the target trial closest to a desired $g^*$.}
As a matter of fact, suppose $A$ and $Z$ are binary with no covariates. Let $g := P(A=1)$, $h := P(Z=1)$, and $b_a := P(Z=1 \mid A=a)$. For any $h^*$ with $h^*_1 = P(Z^*=1)$, the induced marginal is linear in $h^*$: $g(h^*)(1) = g\left(\frac{1-b_1}{1-h}(1-h^*_1) + \frac{b_1}{h} h^*_1\right)$. Thus the attainable set $\mathcal{J}$ is a line segment, and projection amounts to selecting the closest point to $g^*$ on that segment. The operator framework below is the infinite-dimensional analog of this linear map. KL is natural when $g^*$ is treated as a data-generating target and we form a likelihood for pseudo-samples $A^* \sim g^*(\cdot \mid W)$. The $L^2(P_{AW})$ loss instead yields a quadratic risk with a linear operator, enabling a convex optimization viewpoint and gradient-based algorithms.

\subsection{Kullback-Leibler projection}\label{EM}

 {We fix a desired target treatment distribution $g^*(A\mid W)$ and seek the instrument policy $h^\dagger$ whose induced marginal $g(h^\dagger)$ is closest to $g^*$ in the Kullback--Leibler (KL) sense:}
\[
h^\dagger
=
\underset{h^* \in \mathcal{I}}{\arg\min}\,
E_W \left[ \mathrm{KL}\bigl(g^*, g(h^*)\bigr)\right],
\qquad
g(h^*)(A\mid W) = \int p(A\mid Z,W)\,h^*(Z\mid W)\,dZ.
\]

This is equivalent to maximum likelihood estimation (MLE) when $(A_i^*,W_i)$ are drawn with $A_i^*\sim g^*(\cdot\mid W_i)$ and $Z$ is treated as missing data (derivation in the Supplement):
\[
\ell(h^*) = \sum_{i=1}^n \log \int h^*(Z\mid W_i)\,p(A_i^*\mid Z,W_i)\,dZ.
\]

We parameterize $h^*$ using the Highly Adaptive Lasso (HAL) sieve:
\[
\mathcal{I}_n = \left\{
h_\beta(Z=1\mid W) = \frac{1}{1+\exp\!\left(-\sum_{j=1}^{\mathcal{R}_n}\beta_j\phi_j(W)\right)}
: \beta\in\mathbb{R}^{\mathcal{R}_n}, \|\beta\|_1< C_n
\right\}.
\]

\textbf{E-step: }At iteration $k$, form the conditional distribution
\[
p_k(Z\mid A^*,W) =
\frac{\hat p(A^*\mid Z,W)\,h^{(k)}(Z\mid W)}
{\int \hat p(A^*\mid z,W)\,h^{(k)}(z\mid W)\,dz},
\]

\noindent where $\hat p$ is an estimator of the natural density $p(A|Z,W)$. For binary $Z$, define
\[
\tau_i = p_k(Z=1\mid A_i^*,W_i)
= \frac{\hat p(A_i^*\mid 1,W_i)\,h^{(k)}(1\mid W_i)}
{\hat p(A_i^*\mid 0,W_i)\,(1-h^{(k)}(1\mid W_i))+\hat p(A_i^*\mid 1,W_i)\,h^{(k)}(1\mid W_i)}.
\]

The $Q$-function becomes
\[
Q(h\mid h^{(k)})
= \sum_{i=1}^n \Big[ \tau_i \log h(1\mid W_i)
+ (1-\tau_i)\log \{1-h(1\mid W_i)\}\Big].
\]

\textbf{M-step: } Maximize $Q(h\mid h^{(k)})$ over the sieve $\mathcal{I}_n$.  
With HAL/logistic parametrization,
\[
\beta^{(k+1)} = \arg\max_\beta
\sum_{i=1}^n \Big[ \tau_i \log \pi_i(\beta) + (1-\tau_i)\log\{1-\pi_i(\beta)\}\Big]
- \lambda \sum_j |\beta_j|,
\]
where $\pi_i(\beta)=\mathrm{logit}^{-1}\big(\sum_j \beta_j\phi_j(W_i)\big)$ and $\lambda$ is the HAL penalty. Refer to the Supplement for implementation details with \texttt{glmnet} and \texttt{hal9001} \citep{hal9001}.

\begin{algorithm}[h!]
	\caption{EM-HAL algorithm for KL projections}
	\label{alg:em}
	\begin{algorithmic}[1]
		\State \textbf{Input:}  {Target treatment distribution} $g^*(A\mid W)$, density estimate $\hat p(A\mid Z,W)$.
		\State Fix artificial samples $A_i^*\sim g^*(\cdot\mid W_i)$ for $i=1,\dots,n$.
		\State Initialize $h^{(0)}(Z\mid W)$.
		\For{$k=0,1,2,\dots$ until convergence}
		\State \textbf{E-step:} For each $i$, compute
		\[
		\tau_i = p_k(Z=1\mid A_i^*,W_i).
		\]
		\State \textbf{M-step:} Update $h^{(k+1)}$ by solving weighted logistic/HAL regression with soft labels $\tau_i$.
		\EndFor
		\State \textbf{Return:} $h^{(\infty)}$.
	\end{algorithmic}
\end{algorithm}

\subsection{Least squares projection}\label{sec:ls}

We now change the KL divergence in the previous section by the $L^2(P_{AW})$ norm. The finite-dimensional linear map in the binary example extends to a linear operator. Following Proposition \ref{bayes}, define
\[
\begin{aligned}
	B: \quad & \mathcal{F} \longrightarrow \mathcal{G} \\
	& h^*(Z \mid W) \longmapsto [Bh^*](A, W) 
	= g(A \mid W) \, E \left[ \left. \frac{h^*(Z \mid W)}{h(Z \mid W)} \, \right|\, A, W \right],
\end{aligned}
\]
so that for any admissible $h^*$, $Bh^*$ equals the induced marginal $g(h^*)(A \mid W)$. To make this precise, let
\[ \mathcal{F} := \left\{ f \in L^2(P_{ZW}) : \frac{f(Z,W)}{h(Z \mid W)} \in L^2(P_{ZW}) \right\}, \qquad \mathcal{G} := \left\{ g(A\mid W) E\left[ \left. \frac{f(Z,W)}{h(Z\mid W)} \right| A,W \right] : f \in \mathcal{F} \right\}, \]
so that $B$ is well-defined and $\mathcal{G} = \mathcal{R}(B)$ (see Supplement for details). $P_{ZW}$ and $P_{AW}$ denote the joint densities of $(Z,W)$ and $(A,W)$.
Consider the density subsets of $\mathcal{F}$ and $\mathcal{G}$,
\[
\mathcal{F}_1 = \left\{ f \in \mathcal{F} \mid f > 0 \ \mu_Z\text{-a.e.}, \int f \, d\mu_Z = 1 \right\}.
\]
\[
\mathcal{G}_1 = \left\{ f \in \mathcal{G} \mid f > 0 \ \mu_A\text{-a.e.}, \int f \, d\mu_A = 1 \right\}
\]
with $\mu_Z$ and $\mu_A$ common dominating measures.

\begin{propo}\label{inc} $B$ maps densities into densities; i.e. $B(\mathcal{F}_1) \subseteq \mathcal{G}_1$
\end{propo}

We pursue
$$\min_{h^* \in \mathcal{F}_1} \langle g^* -Bh^*,g^* - Bh^* \rangle_{L^2(P_{AW})}$$
(mind the sets over which minimization takes place) so that our risk function is \begin{equation}\label{sq:risk}
	h^* \mapsto R_{g^*}(h^*) = \int (g^* - B(h))^2 dM(A,W)\end{equation}

with $M=P_{AW}$. We arrived to a problem involving minimization of a convex functional over a convex set (find a proof for the next result in the Supplement)

\begin{theo}[Existence and Uniqueness of Minimizer]\label{eu}
	Given the bounded linear operator $B \colon \mathcal{F} \to \mathcal{G} \subseteq L^2(P_{AW})$ and the target function $g^* \in L^2(P_{AW}) $, consider the minimization problem:
	\[
	\min_{h^* \in \mathcal{F}_1} R_{g^*}(h^*), \quad \text{where} \quad R_{g^*}(h^*) = \|g^* - B h^*\|_{L^2(P_{AW})}^2.
	\]
	Assume one of the following holds:
	\begin{enumerate}
		\item $\mathcal{F}_1$ is bounded in $L^2(\mu_Z)$, or
		\item $R_{g^*}$ is coercive on $\mathcal{F}_1$, i.e., $\|h^*\|_{L^2} \to \infty \implies R_{g^*}(h^*) \to +\infty$.
	\end{enumerate}
	Then:
	\begin{enumerate}
		\item[(i)] There exists a minimizer $h^* \in \mathcal{F}_1$ solving the problem.
		\item[(ii)] If $B$ is injective (or $R_{g^*}$ is strictly convex), the minimizer is unique.
		
	\end{enumerate}
\end{theo}

\begin{rem}We note that if $(B^*B)^{-1}B^*g^* \in \mathcal{F}_1$, then that one is the minimizer. 
\end{rem}

Theorem~\ref{eu} states that under certain assumptions, if moreover the operator $B$ is injective then the optimizer $h$ is unique. A natural question is the consequence of $B$ not being injective. Suppose there exist $\tilde h_1, \tilde h_2$ such that $B(\tilde h_1) = B(\tilde h_2)$. Then the implied interventions coincide, $g(\tilde h_1) = g(\tilde h_2)$, and Theorem~\ref{propoim} implies $E[Y^{g(\tilde h_1)}] = E[Y^{g(\tilde h_2)}]$. This is analogous to the recently discussed phenomenology of estimands that are identifiable even though the nuisances needed to compute them are not \citep{neweak}.

\begin{theo}\label{inf}
	The canonical gradient of the theoretical risk \ref{sq:risk} is $D^*(h) = -2 hB^*(g^* - B(h))$
\end{theo}
Note that this is valid for any regular enough linear operator $B$.

Assume the instrument can take $q$ different values so that the probability mass function of $Z$ lives in the $(q-1)-$simplex . Theorem \ref{inf} allows to run a steepest descent algorithm by minimizing risk \ref{sq:risk} along one dimensional paths with score $D^*$ iteratively. After each step, we project the updated $h^*$ onto the positive simplex using the $\mathcal{O}(q \log q)$ (recall $q$ is the dimension of the instrument) procedure by \citep{duchi2008efficient} displayed in Algorithm \ref{alg}, where $q$ is the support size of $Z$.

\begin{algorithm}[h!]
	\caption{Projected canonical gradient descent on the simplex}\label{alg}
	\begin{algorithmic}[1]
		\State \textbf{Input:} Matrix $B$, vector $g^*$, initial point $h^{(0)} \in \mathbb{R}^q$, step size $t > 0$
		\State \textbf{Initialize:} $h^{(0)}$
		\Repeat
		\State Compute canonical gradient at $h^{(k)}$: \[
		D^*(h^{(k)}) = -2 B^\top (g^* - B h^{(k)})
		\]
		\State Gradient step: \[
		h' \gets h^{(k)}- t D^*(h^{(k)})
		\]
		\State Project onto feasible set $\mathcal{F}_1$ using \citep{duchi2008efficient}: \[
		h^{(k+ 1)}\gets \mathcal{P}_{\mathcal{F}_1}(h')
		\]
		\Until{Convergence (e.g., $\|h^{(k+1)} - h^{(k)}\| \leq \epsilon$)}
		
	\end{algorithmic}
\end{algorithm}

\begin{example}
Assume $A,Z$ to be binary and no observed covariates. Considering any $h^*(Z) \in \mathcal{F}_1$ 
	$$ h^*(Z) =\begin{cases}
		h^*_0 & \text{if } Z = 0  \\
		h^*_1  = 1- h^* _0 & \text{if } Z=1
	\end{cases} $$ we have that its image by $B$ is	$B(h^*(Z)) = g(A)\left(\frac{h^*(0)}{h(0)}(1-P(Z=1|A)) + \frac{h^*(1)}{h(1)}P(Z=1|A) \right) $, which amounts to the following matrix times vector multiplication

	\[
	\begin{pmatrix}
		h_0^* \\
		h_1^*
	\end{pmatrix} \overset{B}{\longmapsto}
	\begin{pmatrix}
		\displaystyle\frac{(1 - g)(1 - b_0)}{1 - h} & \displaystyle\frac{(1 - g) b_0}{h} \\
		\displaystyle\frac{g(1 - b_1)}{1 - h} & \displaystyle\frac{g b_1}{h}
	\end{pmatrix}
	\begin{pmatrix}
		h_0^* \\
		h_1^*
	\end{pmatrix}
	\]
	
	with $P(Z^*=1):=h^*_1= 1-h^*_0$, $b_a:=P(Z=1|A=a)$, $g:=g(1)=P(A=1)$, $h=h(1)=P(Z=1)$ and $B$ represented by the matrix.

	 {For illustration, we set the natural policy to $h=0.6$ and $P(A=1\mid Z=z) = 0.3z + 0.5(1-z)$. Empirically, $\hat{b}_0 \simeq 0.68$ and $\hat{b}_1 \simeq 0.47$. For $g^*=(0.65,0.35)$, the unconstrained candidate $(B^*B)^{-1}B^*g^*$, i.e.\ the preimage of the projection of $g^*$ onto $B(\mathcal{F}_1)$, falls inside $\mathcal{F}_1$ and equals approximately $(0.25,0.75)$. However, for $g^*=(0.4,0.6)$ we obtain $(B^*B)^{-1}B^*g^* = (1.49,-0.49) \notin \mathcal{F}_1$, since this quantity is the minimizer of $\langle g^*-Bh^*, g^*-Bh^*\rangle$ over all of $\mathbb{R}^2$ rather than over $\mathcal{F}_1$. Projected gradient descent (Algorithm~\ref{alg}) is therefore required, and yields the constrained minimizer $(1,0)$, with corresponding implied intervention $g(h^*) = (0.5,0.5)$.}

\end{example}

Please refer to the Supplement for estimation of the least squares projection of a desired target treatment distribution with multivariate continuous instruments using score neural networks.

 {\begin{rem}[Inference for the projected estimand]
\label{rem:inference-projection}
The algorithms above yield a data-adaptive instrument policy $\hat{h}^\dagger$, and the downstream causal estimand of interest is $\Psi_{h^\dagger}(P) = E[Y^{h^\dagger}]$, where $h^\dagger$ is itself determined by the data through the projection step. Once $\hat{h}^\dagger$ is fixed, the TMLE of Section~\ref{sec:semi} provides a regular, asymptotically linear estimator of $E[Y^{\hat{h}^\dagger}]$, together with valid confidence intervals for that plug-in quantity. However, the estimand $\Psi_{h^\dagger}(P)$ depends on $P$ also through $h^\dagger$, making it a data-adaptive target parameter in the sense of \citet{hubbard2016statistical}. A full asymptotic treatment---establishing convergence of $\hat{h}^\dagger$ to $h^\dagger$, bounding the resulting second-order remainder, and constructing confidence intervals that account for the variability of the projection step itself---is beyond the scope of the present paper and is left to future work.
\end{rem}}

\section{Numerical examples}

\subsection{G-computation formula}\label{sec:sim1}

We consider a simulation study to evaluate performance of the empirical version of Theorem \ref{idth}, which aims at identifying instrumental stochastic interventions. The data generating mechanism follows the toy NPSEM described in the Supplement. We consider an intervention that modifies the natural distribution of the instrument according to $h^*(Z=1|W=0) = 0.7, h^*(Z=1|W=1) = 0.4$. The target parameter is $\Psi_{h^*}(P) = E[Y^{h^*}]$, the expected outcome under this stochastic intervention on the instrument. The number of replications is $B = 1000$.

We estimate $E[Y|Z,W]$ using ordinary linear least squares regression of $Y$ on $Z$ and $W$. The G-computation estimate is then:
$$\hat{\Psi}^{G}_{h^*} = \frac{1}{n}\sum_{i=1}^n \left[\hat{E}[Y|Z=0,W_i](1-h^*(W_i)) + \hat{E}[Y|Z=1,W_i]h^*(W_i)\right]$$

\begin{figure}[h!]
	\centering
	\includegraphics[width=\textwidth]{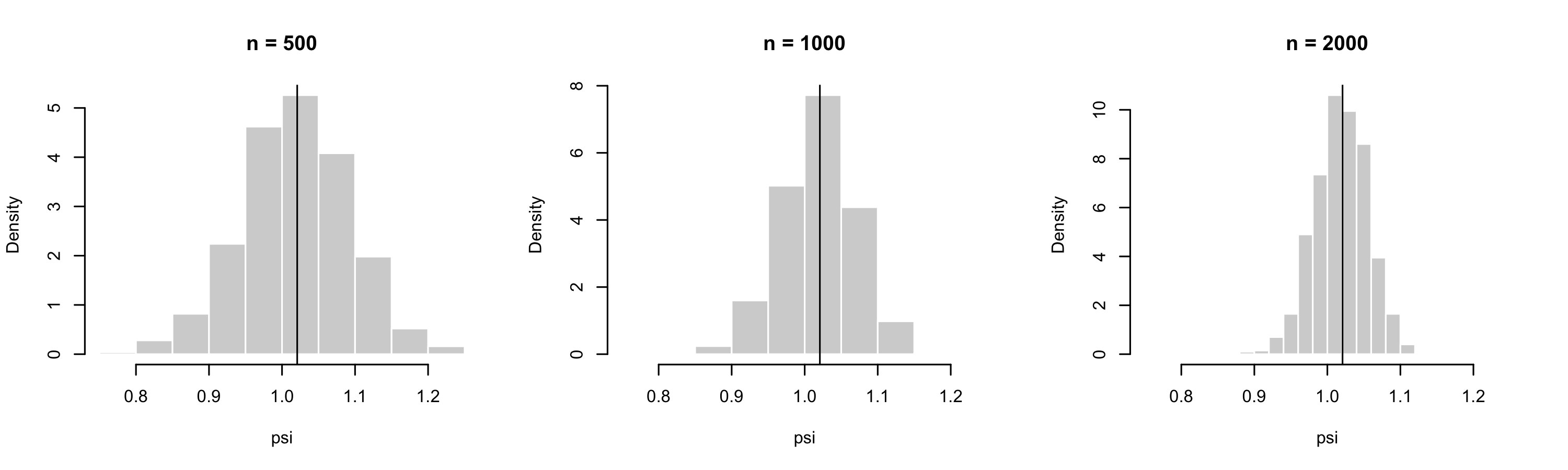}
	\caption{Histograms of estimating $\Psi = EY^{h^*}$ across $1000$ simulation draws for different sample sizes}
	\label{fig:hists}
\end{figure}

The simulation results are visible in Figure \ref{fig:hists} and Table \ref{tab:part1}. 

\begin{table}[H]\centering\begin{tabular}{rrrrrr}  \hline
		$n$ & $\hat \Psi_{\textrm{TMLE}}$ & $\hat \Psi$ & $\hat \sigma$ & Cvrg. $\alpha=0.1$ & Cvrg. $\alpha=0.05$ \\   \hline100 & 1.010 & 0.803 & 0.164 & 0.869 & 0.920 \\   500 & 1.018 & 0.798 & 0.075 & 0.894 & 0.949 \\   1000 & 1.019 & 0.799 & 0.053 & 0.916 & 0.953 \\   2000 & 1.019 & 0.798 & 0.038 & 0.916 & 0.960 \\   10000 & 1.020 & 0.797 & 0.017 & 0.923 & 0.954 \\    \hline\end{tabular}
	\caption{Empirical performance of the G-computation formula in Theorem \ref{idth}. True value is $\Psi = 1.02$}
	\label{tab:part1}\end{table}

\subsection{KL projection}

\begin{figure}[t!]
	\centering
	
	\begin{tabular}{ccc}
		\includegraphics[width=0.3\textwidth]{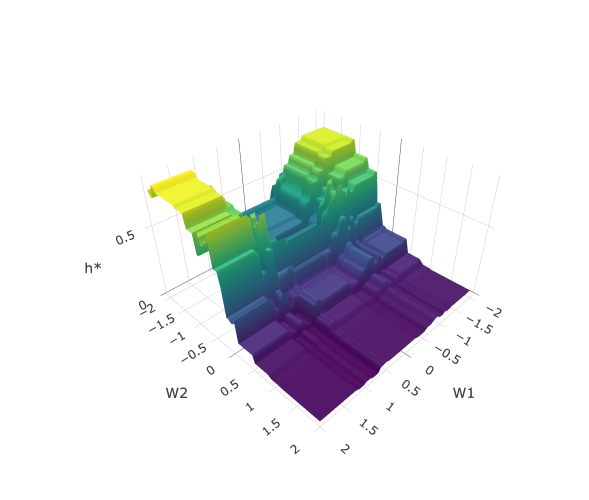} &
		\includegraphics[width=0.3\textwidth]{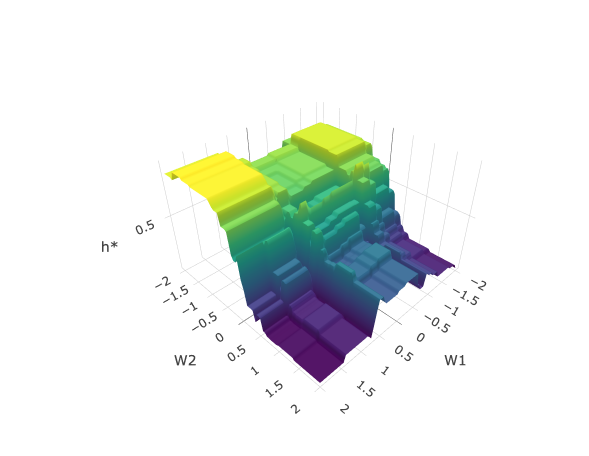} &
		\includegraphics[width=0.3\textwidth]{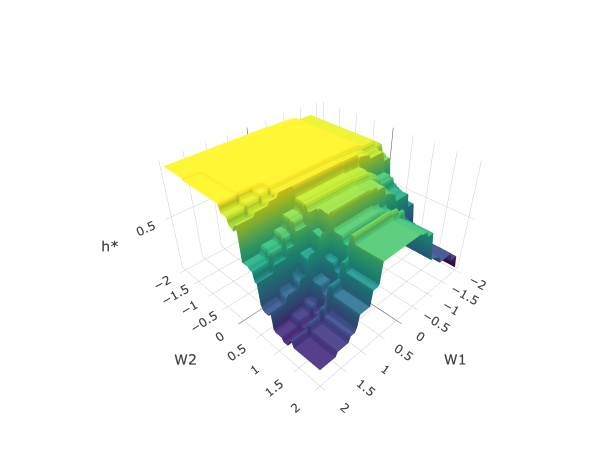} \\
	\end{tabular}

	\begin{tabular}{ccc}
		\includegraphics[width=0.3\textwidth]{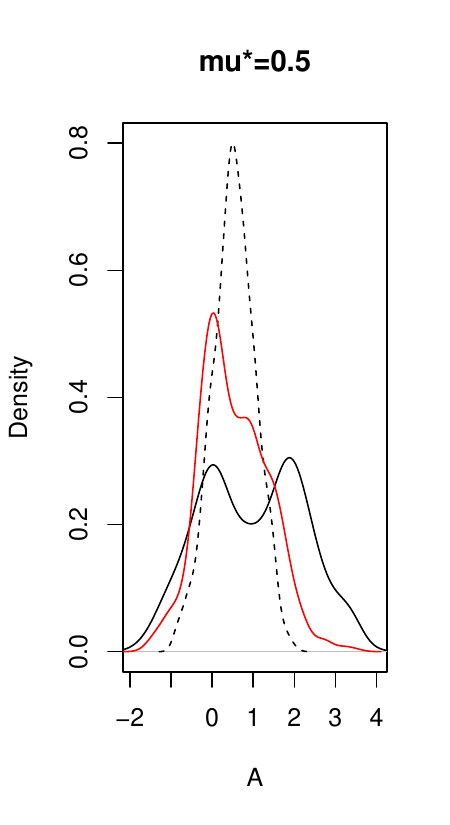} &
		\includegraphics[width=0.3\textwidth]{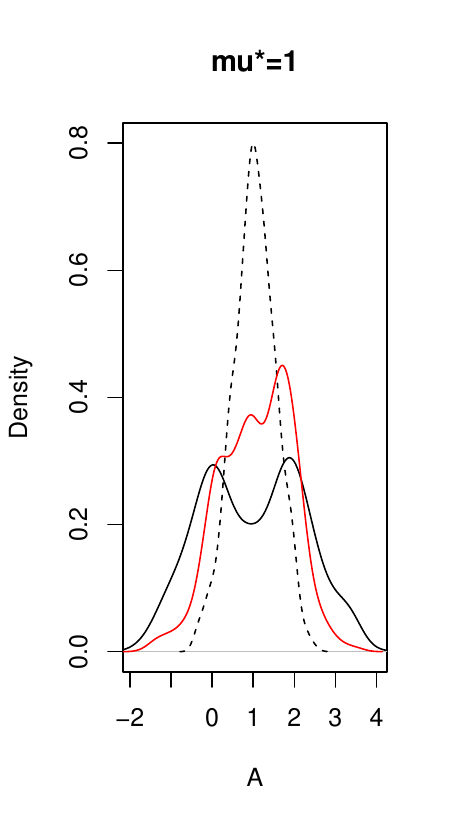} &
		\includegraphics[width=0.3\textwidth]{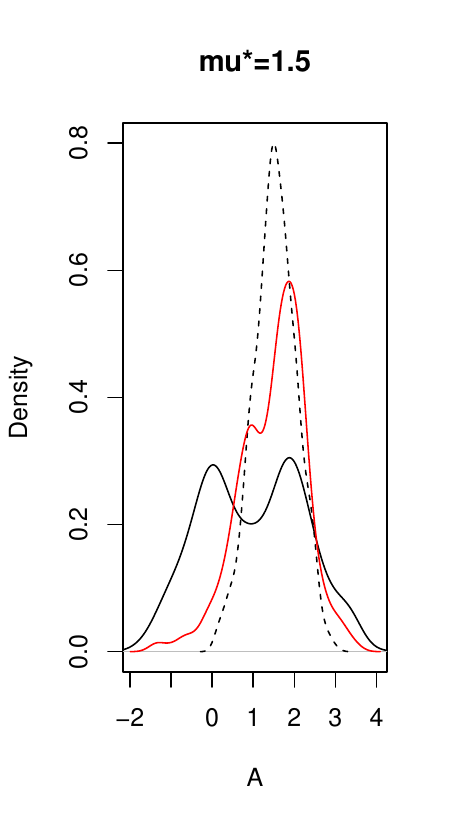} \\
	\end{tabular}
	
	\caption{Simulation results for the EM algorithm, where desired treatment interventions are in dashed, KL-closest implied interventions in red. The surfaces on the first row are the optimal encouragement rules leading to the closest implied intervention, represented as a function of the bivariate covariates. Note how the red line correctly follows the dashed one.}
	\label{fig:combined}
\end{figure}

We implement and test the EM algorithm described in Section \ref{EM} for finding an instrumental intervention $h^*(Z|W)$ that best approximates a desired treatment distribution $g^*(A|W)$ through the relationship established in Theorem \ref{idg}. The covariates are bivariate, with each component uniformly distributed on $[-2,2]$. The instrument $Z \sim \text{Bernoulli}(\text{logit}^{-1}(\phi(W)))$. The target is $g^*(A) = \mathcal{N}(\mu^*, \sigma^{*2})$ with $\mu^* = 1$. The natural propensity of the instrument is $P(Z=1|W) = \text{logit}^{-1}(\phi(W))$ where $\phi(W) = W_1 \sqrt{|W_2|} \text{sign}(W_2)$. The conditional treatment distribution is $A|Z,W \sim \mathcal{N}(\gamma Z + \psi(W), \sigma^2)$ where $\psi(W) = \sin(W_1) \log(1 + W_2^2)$ and $\gamma = 2$. Sample size is $n = 500$ and $\sigma = 0.2$. Figure \ref{fig:combined} compares the implied treatment distribution $A^{h^*}$ (red) with the natural distribution (black) and target (dashed).

\section{Impact of Medicaid on medical visits}\label{sec:oregon}

We use the Oregon Medicaid lottery of 2008 as an instrument to estimate marginal effects of insurance enrollment on primary care utilization. The instrument is lottery selection ($Z=1$ if selected), the treatment is approval status ($A=\texttt{data\$approved\_app}$), and the outcome is any primary care visit within 12 months ($Y=\texttt{doc\_any\_12m}$). Baseline covariates are $\texttt{W.dep\_dx\_0m}$, $\texttt{W.zip\_msa\_list}$, $\texttt{W.birthyear\_list}$, and $\texttt{W.female\_list}$. Sample size is $n=16403$. 

As in Section \ref{hidnohid}, we consider $S'=\texttt{W.dep\_dx\_0m}$.

In Table 2, two distinct interventional stochastic policies of the form $h^*(w) = P(Z^* = 1 \mid W_{S'} = w)$ were specified manually. The G-computation formula from Theorem \ref{idth} was used to estimate $\Psi = E[Y^{h^*}]$, with each component substituted by its corresponding estimate obtained through the TMLE procedure described in Section \ref{sec:semi} including all the 4 relevant baseline covariates in the clever weights and offset, starting from initial estimates fitted using HAL. The reported standard errors $\hat{\sigma}$ are plug-in estimators of $\sqrt{\operatorname{Var}(\text{EIC}) / n}$, where EIC denotes the efficient influence curve derived in Proposition \ref{prop:eic}. 
\vspace{-0.5cm}

\begin{table}[h!]
	\centering
	\begin{tabular}{lcccc}
		\hline
		& $(h^*(0), h^*(1))$ & $\hat{\Psi}$ & $\hat{\sigma}$ & $(g(h^*)(0), g(h^*)(1))$ \\
		\hline
		Policy 1 & $(22.22, 11.11)$ & $58.81$ & $0.45$ & $(7.71, 4.52)$ \\
		Policy 2 & $(77.78, 88.89)$ & $63.92$ & $0.44$ & $(26.98, 36.17)$ \\
		\hline
	\end{tabular}
	\caption{Estimated marginal causal effects $\hat{\Psi}=EY^{h^*}$ under two different stochastic interventions $h^*$, Policy 1 and Policy 2, with estimated standard errors $\hat{\sigma}$ and their implied treatment policies $g^* = g(h^*)$.}
\end{table}

We proceed as in Section \ref{hidnohid}: identification uses a reduced $W$ while TMLE uses all covariates. With $W=\texttt{W.dep\_dx\_0m}$, the empirical kernel is $\widehat P(A = 1 \mid Z = i, W = j) = (0.00, 0.00;35.69,40.69)$ (rows $i$, columns $j$). Thus the null policy $g^*$ with $P(A^{g^*}=1)=0$ is $Z$-compatible, leading to $h^*=0$ by \eqref{hthatimplies}. The natural encouragement design was $(h(0),h(1))=(50.79,46.82)$, while $(h^*(0),h^*(1))=(0,0)$ implies $(g(h^*)(0),g(h^*)(1))=(0,0)$. The results were $\hat{E}Y^0 = 57.49, \hat{\sigma}^0 = 0.53$ and $\hat{E}Y = 61.29, \hat{\sigma} = 0.38$.

We are recovering a counterfactual quantity under almost no assumptions: only proper randomization of the instrument and standard exclusion restriction.  When considering the null policy $h^* = 0$, weak positivity is automatically satisfied.

{
	\section{Discussion}\label{sec:disc}
	
	The central contribution of this paper is to reframe IV identification as a statement about stochastic policies on the instrument and about the {induced} conditional marginals they generate on the treatment. Under the IV core (conditional instrument randomization plus exclusion restriction) and without monotonicity, completeness, or homogeneity, Theorem~\ref{idth} identifies $E[Y^{h^*}]$ for any stochastic intervention $h^*$ on the instrument via a transparent G-computation formula, and Theorems~\ref{propoim}--\ref{idg} show that this effect is mediated through the induced marginal $g(h^*)(A\mid W)$ of the treatment.
	
	An important distinction to be made is the one between \emph{induced marginal} and \emph{free-standing policy}. In a world where one manipulates $Z$ and leaves the rest of the NPSEM untouched, $A$ still depends on $U_A$ and remains correlated with $U_Y$; $g(h^*)$ only describes the treatment distribution that emerges in that world. A separate thought experiment that draws $A$ from $g(h^*)$ independently of $U_Y$ corresponds to a different counterfactual, identified only under additional structure (Remark~S1, Supplementary Material).
	
	 {When a desired target treatment distribution $g^*$ lies outside the range of the induced-marginal map $g(\cdot)$, we project $g^*$ onto the set of attainable induced marginals in KL or $L^2$ sense (Sections~\ref{EM} and~\ref{sec:ls}). These projections yield the identifiable instrument policy whose propagation brings the treatment distribution closest to $g^*$, and the associated causal effect $E[Y^{h^\dagger}]$ is the best identifiable surrogate for the researcher's target. This has a direct interpretation in the target trial emulation framework \citep{hernan2016target}: $\mathcal{J}$ is the set of RCT-like treatment protocols that the IV data can support, and $h^\dagger$ is the instrument policy that emulates the closest feasible trial to the desired $g^*$. Because $\hat h^\dagger$ is selected from the data, inference for the projected estimand $\Psi_{h^\dagger}(P)$ faces the same concern as iterative protocol revision in target trial emulation \citep{stensrud2025target}; this can be addressed through data-adaptive target parameter theory \citep{hubbard2016statistical}, with full development left to future work.}

	 {Instrumental variables give the analyst a lever on the assignment mechanism, not on the treatment itself. When the analyst pulls that lever, a treatment marginal distribution is induced, and this induced distribution---together with the natural confounding structure that remains in the world---defines the only causal quantity that observational data equipped with an instrument can recover without further untestable assumptions. Our framework puts that quantity, and its closest projections to desired treatment-distribution profiles, at the centre of the analysis.}
}

\section*{Data availability statement}

The data that support the findings of this study are openly available from the National Bureau of Economic Research (Oregon Health Insurance Experiment public use data) at \url{https://www.nber.org/research/data/oregon-health-insurance-experiment-data}.


\section{Proof of Theorem 1}

$$\begin{aligned}
	EY^{h^*} = E( E(Y^{h^*} |W)) =  E( 	E(E[Y^{h^*}|Z^*,W] |W)) 
\end{aligned}$$
Define  $$\varphi(Z^*,W):=E[f_Y(f_A(Z^*,W,U_A),W,Z^*,U_Y)|Z^*,W]$$ 
Then we have

$$\begin{aligned}
	E(E[Y^{h^*}|Z^*,W] |W) = E_{Z^*|W \sim h^*(W)}\varphi(Z^*,W)
\end{aligned}$$

so that the only randomness left is due to $U_A$ and $U_Y$, which are conditionally independent of $Z^*$ given $W$ because $Z^*$ was properly randomized. Therefore,

$$\varphi(Z^*,W)=E[f_Y(f_A(z^*,W,U_A),W,z^*,U_Y)|W] \Big|_{z^*=Z^*}$$

On the other hand
$$E(Y|Z,W)=E[f_Y(f_A(Z,W,U_A),W,Z,U_Y)|Z,W] = E[f_Y(f_A(z,W,U_A),W,z,U_Y)|W] \Big|_{z=Z}$$

where the last equality holds because we are assuming that the natural law $P$ is such that $Z$ is properly randomized ($Z$ independent of residuals $U_A,U_Y$ given $W$)

So we have $E(Y|Z,W) = \varphi(Z,W)$ and the final G computation formula is 

$$EY^{h^*} = E_WE_{(Z^*|W) \sim h^*(W)}E(Y|Z=Z^*,W)$$

\section{Proof of Proposition 4}

$$\begin{aligned}\int  g(A|W){E} \left[\left. \frac{h^*(Z|W)}{h(Z|W)} \,\right|\, A, W \right] dA = \int  g(A|W)\int  \frac{h^*(Z|W)}{h(Z|W)} p(Z|A,W)dZdA = \\ 
	= \int \int \frac{p(A,W)p(W)p(Z,A,W)h^*(Z|W)}{p(W)p(A,W)p(Z,W)} dZdA =\int \int \frac{p(Z,A,W)h^*(Z|W)}{p(Z,W)} dZdA = \\ \int \int {p(A|Z,W)h^*(Z|W)}dZdA = \int h^*(Z|W)\left(\int {p(A|Z,W)}dA\right)dZ = 1  \end{aligned}$$

\section{Toy example from the Introduction}

We illustrate the induced-marginal map with a simple binary NPSEM:
\[
W = 1 \{U_W < 0.3\},\quad
Z = 1 \{U_Z < h(W)\},\quad
A = 1 \{U < p(Z,W)\},\quad
Y = 2A + W - U + U_Y,
\]
with $h(0)=0.3$, $h(1)=0.8$, and $P(A=1\mid Z=i, W=j):=p(i,j)$ given by the $ij$th entry of the $2 \times 2$ matrix $(0.3,0.8,0.7,0.5)$ (row-major). The errors $U_W, U_Z, U \sim \mathrm{Uniform}(0,1)$ and $U_Y \sim \mathcal{N}(0,0.05^2)$ are mutually independent, and $U_Z$ is independent of $(U,U_Y)$ given $W$.

For an intervention on the instrument, $h^*(W) = P(Z^{h^*}=1\mid W)$, the induced marginal is
\[ g(h^*)(A=1\mid W=w) = p(0,w)\{1-h^*(w)\} + p(1,w)h^*(w). \]
Two intervened NPSEMs used for validation are
\begin{equation}
\begin{aligned}
W&  = 1 \{U_W < 0.3 \}, \\
Z^{h^*}&= 1 \{U_Z < h^*(W)\}, \\
A^{h^*} &= 1 \{U < p(Z^{h^*},W)\}, \\
Y^{h^*} &= 2A^{h^*} + W - U + U_Y,
\end{aligned}
\end{equation}
and
\begin{equation}
\begin{aligned}
W& = 1 \{U_W < 0.3 \}, \\
A^{g^*} &= 1 \{U < g^*(W)\}, \\
Y^{g^*} &= 2A^{g^*}  + W - U + U_Y.
\end{aligned}
\end{equation}

For the choice $h^*(0)=0.7$, $h^*(1)=0.4$, the population G-computation formula yields $E[Y^{h^*}] = 1.02$. In simulations with $10^6$ samples, the empirical averages are $1.0203$ for $Y^{h^*}$ and $1.0189$ for $Y^{g^*}$, matching Theorem 2. A final simulation run from the natural NPSEM gives a marginal mean of $Y$ of $0.7248$ with $P(A=1\mid W=w) = (0.42, 0.56)$ for $w=0,1$.

All numerical experiments and algorithmic implementations described in this work can be reproduced using the code at \url{https://github.com/meixide/implied_interventions}.

\section{Induced Marginal vs.\ Free-Standing Policy on \texorpdfstring{$A$}{A}}

\begin{rem}[Induced marginal vs.\ free-standing policy on $A$]\label{rem:indep}
Let $\widetilde Y^{g^*} = f_Y(\widetilde A^{g^*},W,U_Y)$ with $\widetilde A^{g^*}\sim g^*(\cdot\mid W)$ drawn \emph{independently} of $U_Y$ given $W$. Then $E[Y^{h^*}] = E[\widetilde Y^{g(h^*)}]$ would hold under either of the following conditions:
    \begin{enumerate}
        \item[(i)] {No hidden confounding:} $U_A \perp U_Y \mid W$. In this case $A^{h^*}$ and $\widetilde A^{g(h^*)}$ have the same joint law with $U_Y$ given $W$.
        \item[(ii)] {Additive outcome:} $f_Y(A,W,U_Y) = m(A,W) + v(W,U_Y)$. Then
        \[
        E\!\left[Y^{h^*}\right] \;=\; E_W\!\left\{\int m(a,W)\, g(h^*)(a\mid W)\, d\mu_A(a) + E[v(W,U_Y)\mid W]\right\} \;=\; E\!\left[\widetilde Y^{g(h^*)}\right],
        \]
        because the $U_Y$-channel is additively separable from $A$ and therefore insensitive to whether $A$ is correlated with $U_Y$ through $U_A$.
    \end{enumerate}
    Our identification result targets $E[Y^{h^*}]$, the intention-to-treat-style effect of a stochastic encouragement, and this is what we claim throughout. The free-standing policy reading is a corollary available only when extra structure holds.
\end{rem}

\section{Proof of Theorem 5}

\begin{lemma}\label{generald}
	Let $Z$ and $A$ be binary and $g^*$ be any treatment intervention satisfying weak monotonicity. Under the hypothesis of Theorem 1, 
	
	$$EY^{g^*} = E\left[ g^*(W) \operatorname{Wald}(W) \right]+ E \left[ \frac{E[Y| Z=0,W]g_1(W)- E[Y|Z=1,W]g_0(W)}{g_1(W)-g_0(W)}\right]$$

\end{lemma}

\begin{proof}
$$\begin{aligned}&EY^{g^*} = \int_W\frac{E(Y|Z=0,W)(g_1(W) - g^*(W)) + E(Y|Z=1,W)(g^*(W) - g_0(W))}{E(A |Z=1,W) - E(A |Z=0,W)} dP(W) = \\
	&= \int_Wg^*(W)\frac{E(Y|Z=1,W)- E(Y|Z=0,W)}{E(A |Z=1,W) - E(A |Z=0,W)} + \\&+ \int_W\frac{E(Y|Z=0,W)g_1(W)- E(Y|Z=1,W)g_0(W)}{E(A |Z=1,W) - E(A |Z=0,W)}
	dP(W)\end{aligned}$$
	\end{proof}
Theorem 5 is proven by observing that the second term on the RHS of the result in Lemma 1 does not depend on $g^*$.

\section{Proof of Theorem 7}

{(i) Existence:} We verify the conditions of Proposition 1.2 from \\
\cite{Ekeland1999}: {Convexity:} $R_{g^*}$ is convex because the squared $L^2$-norm $h^* \mapsto \|g^* - B h^*\|^2$ is convex and $B$ is linear, so $h^* \mapsto B h^*$ is affine. {Lower Semicontinuity (lsc):} The $L^2$-norm is weakly lsc, and $B$ is bounded (hence continuous), so $R_{g^*}$ is lsc with respect to weak convergence in $L^2(\mu)$. Feasibility: $\mathcal{F}_1$ is nonempty and closed under weak convergence in $L^2(\mu_Z)$ by Fatou's lemma. If $\mathcal{F}_1$ is bounded, its weak closure is weakly compact (by Banach-Alaoglu in reflexive spaces). If $R_{g^*}$ is coercive, sublevel sets $\{h^* \mid R_{g^*}(h^*) \leq \alpha\}$ are bounded and weakly closed. Hence, by Proposition 1.2 in \cite{Ekeland1999}, a minimizer exists.
(ii) Uniqueness: If $B$ is injective, the map $h^* \mapsto B h^*$ is strictly convex on $\mathcal{F}_1$, making $R_{g^*}$ strictly convex. Uniqueness follows from the fact that strictly convex functions have at most one minimizer on convex sets. 

\section{Proof of Theorem 8}
Let $\tilde h = \tilde h (g^*)$ be the minimizer of the constrained problem. We know $\tilde h$ verifies $$\frac{d}{d\delta} R_{g^*}((1+\delta S) \tilde h) |_{\delta=0} = 0$$ for all $S$ verifying $E_{\tilde h}(S(Z|W) |W)=0$. Now we define a perturbation path given by an arbitrary score $S$:
$h_\delta = (1 + \delta S) h$
and differentiate the risk:

$$\left. \frac{d}{d\delta} R_{g}(h_\delta) \right|_{\delta=0}
= \left. \frac{d}{d\delta}{E}[(g^* - B(h_\delta))^2] \right|_{\delta=0}  = -2 \, {E}[(g^* - B(h)) B(Sh)]$$
Therefore the pathwise derivative is:
$$ -2 \langle g^* - B(h), B(S h) \rangle_{L^2(P_{AW})}$$
Working this expression out

$$\begin{aligned}-2 \langle g^* - B(h), B(S h) \rangle_{L^2(P_{AW})}  =  -2 \langle B^*(g^* - B(h)), S h \rangle_{L^2(P_{AW})}  \\= -2 \langle  hB^*(g^* - B(h)), S \rangle_{L^2(P_{AW})} = E_{AW}(D^*S)\end{aligned} $$

\begin{rem}
	Changing the measure $M$ is also going to change the projection, as we are changing the norm. The natural choice is $M(A,W)=P(A,W)$ the joint density of treatment and covariates. Another option could be the choice $dM =\frac{dP(W)}{g(A|W)}dA$, because if we pull out $g^2$ the risk function becomes $$E_{P_{AW}}\left(E\left[\frac{h^*}{h} \middle| A,W\right]-\frac{g^*}{g}\left(A|W\right) \right)^2$$ involving both Riesz representers for $EY^{h^*}$ and $EY^{g^*}$. 
\end{rem}

\section{Proof of Proposition \ref{clashift}}

We focus our attention in interventions of the form \cite{diaz} 

$$h^*_\beta(Z|W) = h(Z - \beta(W)|W)$$

with $\beta(W)$ having the same dimension as $Z$

Assume $h(Z \mid W)$ is differentiable. Then, for small $\beta(W)$, we can approximate the shift via a Taylor expansion of the log-density:
\[
\log h(Z - \beta(W) \mid W) = \log h(Z \mid W) - \beta(W) ^T \psi(Z,W) + O(\beta(W)^2),
\]
with
\[
\psi(Z,W) = \frac{\partial}{\partial Z} \log h(Z \mid W).
\]
which is the score neural network of a diffusion model \cite{feng2024optimalconvexmestimationscore}
Exponentiating both sides yields:
\[
h(Z - \beta(W) \mid W) \approx h(Z \mid W) \cdot \exp\left( - \beta(W) ^T \psi(Z,W) \right).
\]

This resembles an exponential tilt, with:
\[
\varphi(Z, W) = - \beta(W) ^T  \psi(Z,W).
\]

Thus, the shifted density can be approximated as:
\[
h^*(Z \mid W) \propto h(Z \mid W) \cdot \exp\left( -\beta(W)^T\psi(Z,W) \right).
\]

Using the exponential tilting representation, the density ratio becomes:
\[
\frac{h^*(Z \mid W)}{h(Z \mid W)} =C(W){\exp\left(-\beta(W)^T \psi(Z,W)\right)}
\]

\section{Score matching through diffusion models}\label{app:diff}

The score matching framework  introduced by \cite{jmlrscore} provides an alternative to maximum likelihood estimation by targeting the pdf derivative instead of the pdf itself. For a conditional pdf $h_0(z|w)$ on $\mathbb{R}^q$, the conditional score function is defined as the gradient of the log-density with respect to the observations:

\begin{equation}
	\psi_0(z,w) := \nabla_z \log h_0(z|w) = \frac{\nabla_z h_0(z|w)}{h_0(z|w)} {1}{\{h_0(z|w) > 0\}}.
\end{equation}

The population score matching objective minimizes the expected squared Euclidean $q$-norm between a candidate score function $\psi$ and the true score $\psi_0$:

\begin{equation}
	\psi \longmapsto	{E}_{h_0}\big(\|\psi(z,w) - \psi_0(z,w)\|^2 \big).
\end{equation}

Through integration by parts, this can be reformulated as minimizing the score matching objective:

\begin{equation}
	\psi \longmapsto {E}_{h_0} \big( \|\psi(z,w)\|^2 + 2 (\nabla_z \cdot \psi)(z,w) \big),
\end{equation}

where $\nabla_z \cdot \psi$ denotes the divergence with respect to $z$. This formulation circumvents explicit computation of the density. 

The implementation employs denoising score matching, where the network learns to estimate the score function of a perturbed version of the data. For Gaussian noise $\epsilon \sim \mathcal{N}(0,\sigma^2 I)$, the target score becomes:

\begin{equation}
	\psi_{\text{target}}(z,w) = -\frac{\epsilon}{\sigma^2} = -\frac{z - z_{\text{clean}}}{\sigma^2}.
\end{equation}

The neural network architecture takes as input the concatenated pair $(z,w)$ and outputs an estimate of $\nabla_z \log h(z|w)$. The training procedure minimizes the mean squared error between the network's predictions and the target scores:

\begin{equation}
	\theta \longmapsto {E}\big[\|\psi_\theta(z+\epsilon,w) - \psi_{\text{target}}(z+\epsilon,w)\|^2\big].
\end{equation}

where $\theta$ belongs to a neural network class.

This approach demonstrates several advantages highlighted in the theoretical framework: it handles non-log-concave distributions effectively, provides robust estimation for heavy-tailed error distributions, and can achieve favorable asymptotic efficiency compared to traditional convex $M$-estimators. 

\begin{figure}[h!]
	\centering
	\includegraphics[width=\textwidth]{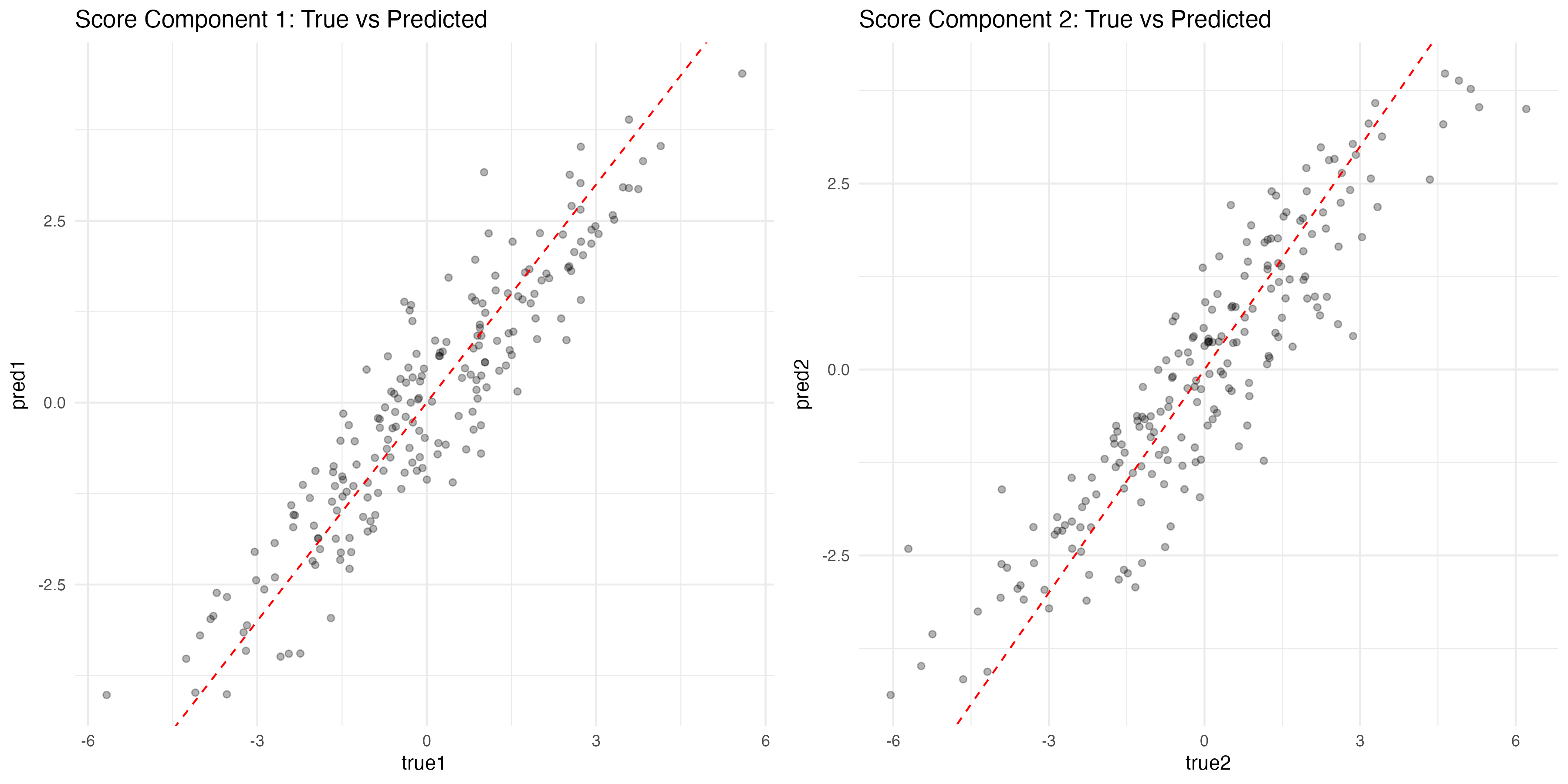}
	\caption{Estimates of the scores $\frac{\partial}{\partial Z} \log h(Z \mid W)$ vs true ones using a diffusion model. Note that $n=200$ only. }
	\label{fig:scores}
\end{figure}

\section{Extended literature review}

  The second, more classical strand of the IV literature emerges from the econometric tradition of simultaneous equations, where the instrument helps resolve endogeneity arising from equilibrium relationships among economic variables. According to \citep{imbens2014instrumental}. ``\textit{simultaneous equations models are both at the core of the econometrics canon and at the core of the confusion concerning instrumental variables methods in the statistics literature}''.
  Although non additive versions can also be found in the literature \citep{matzkin2003nonparametric, imbens2009identification}, this line of work often assumes $f_Y(A,W,U_Y) = f_P(A,W) + U_Y$, where $f_P$ solves the conditional moment equation $E [Y | Z]= E [f(A,W) | Z]$ derived from exclusion restriction $E [\varepsilon_Y | Z] = 0$. 
    By definition, identifiability of $f$- which means being able to write it as a functional of the natural probability law $P$- takes place when the previous Fredholm type I integral equation is such that a solution for it exists and is unique. By direct inspection, uniqueness is translated into being the dependence of $A$ and $W$ on $Z$ strong enough so that different $f$s are not sent to similar functions of $Z$ after applying the conditional expectation operator $ f \longmapsto E [f(A,W) | Z]$. \citep{newey03} formalized this characterization of uniqueness through  completeness \citep{lehmann1950} of the endogenous regressors (including the treatment) conditional pdf given the instrument, which can be easily violated when the IV are weak. Also, \citep{canay2013testability} point out that completeness is impossible to be tested under general nonparametric models. In Banach spaces language, the solution to the previous equation is unique if and only if the null space of the conditional expectation operator is zero. Introducing assumptions that turn this operator and its adjoint into Hilbert-Schmidt operators, \citep{darolles} claim that $f$ is identifiable if and only if the spectrum of the associated information operator does not contain zero (equiv. being one-to-one). 
    
    This equivalence holds only under an important caveat: namely, that a solution to the inverse problem exists—equivalently, that the conditional expectation ${E}(Y \mid Z)$ lies in the range of the conditional expectation operator. Independently of this, \citep{darolles} show that a solution exists if and only if a certain series involving the reciprocals of the eigenvalues $\lambda_i$ of the operator converges. However, since the operator is compact, its eigenvalues tend to zero, and the presence of $1/\lambda_i$ in the theoretical solution renders the inverse problem ill-posed. Despite this, the authors construct a kernel-based estimator $f_n^{\alpha_n}$ and prove consistency of $\| f_n^{\alpha_n} - f \|$, contingent on specific rates of convergence involving the regularization parameter $\alpha_n$ and bandwidths. Crucially, this requires that the bias from regularization—referred to as the identifiability bias—vanishes as $\alpha_n \to 0$ (Assumption A.2 in \citep{darolles}). Yet, this assumption is difficult to verify or interpret in practice, leading to conceptual tension: although the existence of a solution $f_P$ may be theoretically questionable, it nonetheless serves as the limiting object of a consistent estimator under stringent and opaque conditions. \citep{ghassami2025debiasedillposedregression} propose a debiased estimator for $f_P$ based on influence function debiasing of a projected mean squared error objective. Their estimator improves robustness to misspecification of nuisance components and achieves favorable finite-sample convergence rates under mild assumptions. 
    
    In contrast, \citep{neweak} take a fundamentally different perspective on inference under weak identification. Rather than aiming to recover the ill-posed nuisance function itself (such as the NPIV regression function), they focus on inference for linear functionals—for example, the average treatment effect—defined over these functions. Their key contribution is to identify conditions under which the functional is strongly identified even if the underlying function is not. In particular, they establish asymptotic normality at $\sqrt{n}$-rates for debiased estimators of such functionals, using a penalized minimax procedure to estimate both the primary and so-called “de-biasing” nuisance functions. However, this methodology still relies on verifying certain high-level regularity conditions and constructing valid debiasing functions that may require knowledge about the operator and data-generating process that is difficult to ascertain in practice.  {\citep{vanderlaan2025nonparametricinstrumentalvariableinference} introduce  projection parameters for linear inverse problems that remain valid without identification, establishing asymptotic results for linear functionals with many weak instruments and proposing cross-fold debiased estimators.}

\section{Technical details for Kullback-Leibler projection}\label{EM-app}


We start with the following definition of closest implied intervention to a given $g^*$ as a projection on $\mathcal{J}$ with respect to KL divergence, which corresponds to minimizing over $h^* \in \mathcal{I}$ 
$$ \operatorname{K L}\left(g^*,g(h^*)\right) = \int g^*(A|W) \log \frac{g^*(A|W)}{g(h^*)(A|W)}dA = E_{g^*(\cdot |W)}\left[\log \frac{g^*(A|W)}{g(h^*)(A|W)} \middle| W\right]$$

where $g$ is the map defined in Theorem 3 so that the implied intervention of the minimizer will be the closest implied intervention to $g^*$ in the KL sense. Working out the previous expression, 
$$\begin{aligned}
	 \operatorname{K L}\left(g^*,g(h^*)\right)  =  \int g^*(A|W) \log{g^*(A|W)}dA -  \int g^*(A|W) \log{g(h^*)(A|W)}dA
	\end{aligned}$$
	
	We can ignore the first term for optimization purposes as it does not involve $h^*$, so that 
		$$\begin{aligned}
	\underset{h^* \in \mathcal{I}}{\arg \min} \,	\operatorname{K L}\left(g^*,g(h^*)\right)  &=  \underset{h^* \in \mathcal{I}}{\arg \max} \, \int g^*(A|W) \log{g(h^*)(A|W)}dA = \\ &=\underset{h^* \in \mathcal{I}}{\arg \max} \, \int g^*(A|W) \left[\log{\int p(A|Z,W) h^*(Z|W)}dZ\right]dA \\ &=\underset{h^* \in \mathcal{I}}{\arg \max} \, E_{g^*(\cdot |W)} \left[\log{\int p(A|Z,W) h^*(Z|W)}dZ \middle| W\right] 
	\end{aligned}$$

For all observed $W_i \in \mathbb{R}^d$, an artificial observation $A^*_i$ can be drawn from the desired $g^*(\cdot|W_i)$. The full data is now $(W,Z^*,A^*)$, but we only observe $(W,A^*)$ where $A^*|W$ is distributed according to $g^* $. Plugging in the empirical $\frac{1}{n} \sum_i\delta_{(A^*_i,W_i)}$ to replace $g^*$ in the last expression shows that the KL-divergence projection of the NPMLE over the class $\mathcal{J}$ coincides with the Maximum Likelihood estimator, 
$$  \underset{h^* \in \mathcal{I}}{\arg \max} \, \sum_{i=1}^n\log \left(\int h^*(z|W_i) p(A^*_i|z,W_i) dz\right) $$

where $(A^*_i,W_i)$ is artificially sampled as $(A^*_i\sim g^*(\cdot,W_i)$. Finally, we need to write down a model for $h^*$, i.e. specifying what the index class $\mathcal{I}$ is. We choose the Highly Adaptive Lasso (HAL) sieve \citep{hal}.  \begin{equation}\label{sieve}
	\mathcal{I}_n = \left\{C(\beta)\exp{\left(\sum_{j} \beta(j) \phi_j\right) : \beta(j) \in \mathbb{R}, j = 1,\ldots,\mathcal{R}_n }\right\}\end{equation}
	so that MLE becomes optimization over $\beta(j), 1 \leq j \leq \mathcal{R}_n $. 
 If $Z$ is binary, then the basis functions $\phi_j \equiv \phi_j(W)$ just depend on $W$ and therefore $\mathcal{R}_n \leq n(2^d-1)$. $C(\beta)$ is a normalizing constant ensuring that the elements in $\mathcal{I}_n$ integrate to 1.

The underlying full data likelihood is modelled as $p(A|Z,W)h^*(Z|W)dP(W)$ where $p(A|Z,W)$ is already given. Therefore, our full data density is parametrized by $h^*(Z|W)$. We know how to perform full data MLE based on a full-data sample $(Z^*,W)$, so that this makes us propose an Expectation-Maximization (EM) routine \citep{dempster}.

The observed log-likelihood for one datum at EM step $k$ evaluated on a generic $h^*\in \mathcal{I}_n$ is
$$\begin{aligned}
	\log \int p(A^*|z,W) h^*(z|W)dz   &=  \log \int \frac{ p(A^*|z,W) h^*(z|W)p_k(z|A^*,W)}{p_k(z|A^*,W)}dz   \\  &= \log E_{p_k} \left[\frac{ p(A^*|Z^*,W) h^*(Z^*|W)}{p_k(Z^*|A^*,W)} \middle| A^*,W\right] \\ &\geq  E_{p_k} \left[\log \frac{ p(A^*|Z^*,W) h^*(Z^*|W)}{p_k(Z^*|A^*,W)} \middle| A^*,W\right] \\ &=   E_{p_k} \left[\log \frac{ p(A^*|Z^*,W) h^*(Z^*|W)}{p_k(Z^*|A^*,W)} \middle| A^*,W\right]. \\
\end{aligned}$$

\noindent Therefore, for E step it suffices to compute 
$$ \sum_{i=1}^n  E_{p_k} \left[\log {h^*(Z^*|W)} \middle| A^*,W\right] \Big|_{A^*_i,W_i}$$

\noindent as maximization of this functional takes place on $h^* \in \mathcal{I}_n$ and we can ignore the terms not depending on $h^*$, which are additive after taking the log of products and divisions. The latter expectation is taken with respect to
$$  p_k(Z^*|A^*,W) = \frac{p_k(A^*,Z^*,W)}{p_k(A^*,W)} =  \frac{p_k(A^*|Z^*,W) p_k(Z^*,W)}{p_k(A^*,W)} =  \frac{p(A^*|Z^*,W) h_k(Z^*|W)}{p_k(A^*|W)}$$

\noindent where $p_k=  g(h_k)$, so that what comes next is pure MLE, leading to M-step 

$$ \underset{h^* \in \mathcal{I}_n}{\arg \max } \sum_{i=1}^n  E_{h_k} \left[   \frac{p(A_i^*|Z^*,W_i) }{g(h_k)(A^*_i|W_i)} \log {h^*(Z^*|W)}\middle| W=W_i\right]$$

\noindent where the $A^*_i$ inside of the expectation is not integrated out.

EM involves repeatedly carrying out E and M steps. E step consists in estimating the conditonal expectation of the full data log-likelihood given the observed data, taken under the full data current estimate. Then this conditional expectation can be replaced by an empirical mean over draws of $Z^*$, becoming exactly like a full data MLE. Therefore, the M step is the same as computing a full data MLE based on the imputed data set and using weights for each imputed observation that all add up to $1/n$. 

\subsection{Implementing the M-step with \texttt{glmnet}}

The EM M-step requires maximization of the $Q$-function
\[
Q(\beta) \;=\; \sum_{i=1}^n \Big\{ \tau_i \,\log \pi_i(\beta) 
+ (1-\tau_i)\,\log\bigl(1-\pi_i(\beta)\bigr) \Big\}, 
\qquad
\pi_i(\beta)=\mathrm{logit}^{-1}\!\big(\beta^\top \phi(W_i)\big),
\]
where $\tau_i = P_k(Z=1\mid A_i^*,W_i)$ are the E-step posteriors.  
This is the log-likelihood of a logistic regression problem with \emph{fractional responses} $\tau_i$.  

However, standard software such as \texttt{glmnet} does not allow fractional labels directly.  
It expects either binary responses $y\in\{0,1\}$ or grouped binomial counts.  
A convenient workaround is the \emph{duplication trick}: for each observation $i$ we create two pseudo-observations,
\[
(y_i=1,\;\text{weight}=\tau_i), 
\qquad
(y_i=0,\;\text{weight}=1-\tau_i).
\]
The weighted log-likelihood contribution of these two rows is
\[
\tau_i\log \pi_i(\beta) + (1-\tau_i)\log(1-\pi_i(\beta)),
\]
which is exactly the $i$th summand of $Q(\beta)$.  
Thus duplicating the data recovers the correct M-step objective.  

\medskip
In \texttt{R}, given a design matrix $X$ of HAL basis functions $\phi_j(W_i)$, 
a vector of soft labels $\tau$, and weights $w=(\tau,1-\tau)$, one can call
\begin{verbatim}
	library(glmnet)
	
	# Expand dataset
	X2 <- rbind(X, X)
	y2 <- c(rep(1, nrow(X)), rep(0, nrow(X)))
	w2 <- c(tau, 1 - tau)
	

	fit <- glmnet(
	x = X2, y = y2,
	family = "binomial",
	weights = w2,
	alpha = 1,          
	standardize = FALSE,  
	intercept = FALSE    
	)
\end{verbatim}
Here the coefficient vector $\hat\beta$ solves the penalized M-step problem
with $\ell_1$ penalty.  Fixing a value of \texttt{lambda} across EM iterations 
preserves the ascent property of the algorithm.

\section{Domain and codomain of \texorpdfstring{$B$}{B}}\label{app:fg}
Define
$$\mathcal{F}:=\left\{f \in L^2\left(P_{Z W}\right): \frac{f(Z, W)}{h(Z \mid W)} \in L^2\left(P_{Z W}\right)\right\} \subseteq L^2(P_{ZW})$$
Using Jensen's inequality we have that

$${E}\left[\left({E}\left[\left.\frac{f(Z, W)}{h(Z \mid W)} \right\rvert\, A, W\right]\right)^2\right] \leq {E}\left[{E}\left[\left.\left(\frac{f(Z, W)}{h(Z \mid W)}\right)^2 \right\rvert\, A, W\right]\right] ={E}\left[\left(\frac{f(Z, W)}{h(Z \mid W)}\right)^2\right]<\infty .$$
Therefore $$f \in \mathcal{F} \Longrightarrow \left\|{E}\left[\left(\frac{f(Z, W)}{h(Z \mid W)}\right)^2\right]\right\|_{L^2} \leq\left\|\frac{f(Z, W)}{h(Z \mid W)}\right\|_{L^2}$$
so that it makes sense to define the following set
$$\mathcal{G}=\left\{g(A|W){E}\left[\left.\frac{f(Z, W)}{h(Z \mid W)} \right\rvert\, A, W\right]: f \in \mathcal{F}\right\} \subseteq L^2\left(P_{A W}\right) $$
where the inclusion holds assuming $g(A|W)$ is bounded. Consider the standard $L^2(P_{AW})$ inner product restricted to $\mathcal{G}$
$$\left\langle g_1, g_2\right\rangle_{L^2(P_{AW})}:={E}_{A W}\left[g_1(A, W) g_2(A, W)\right]=\int g_1(a, w) g_2(a, w) d P_{A W}(a, w)$$

$\mathcal{G}$ is a closed subspace of $L^2\left(P_{A W}\right)$ because the conditional expectation ${E}[\cdot \mid A, W]$ is an orthogonal projection.
Thus, $\mathcal{G}$ inherits completeness from $L^2\left(P_{A W}\right)$. Therefore, 
\begin{propo} $\left( \mathcal{G}, \left\langle \cdot,\cdot\right\rangle_{L^2(P_{AW})}\right)$ 
	is a Hilbert space
\end{propo}

Consider the density subsets.
	\[
\mathcal{F}_1 = \left\{ f \in \mathcal{F} \mid f > 0 \ \mu_Z\text{-a.e.}, \int f \, d\mu_Z = 1 \right\}.
\]
\[
\mathcal{G}_1 = \left\{ f \in \mathcal{G} \mid f > 0 \ \mu_A\text{-a.e.}, \int f \, d\mu_A = 1 \right\}
\]
with $\mu_Z$ and $\mu_A$ common dominating measures.

\begin{propo}\label{inc-app} $B$ maps densities into densities; i.e.
	$$B(\mathcal{F}_1) \subseteq \mathcal{G}_1$$
	\end{propo}

\begin{example}	Assume $A,Z$ to be binary and no observed covariates. Considering any $h^*(Z) \in \mathcal{F}_1$ 
	$$ h^*(Z) =\begin{cases}
		h^*_0 & \text{if } Z = 0  \\
		h^*_1  = 1- h^* _0 & \text{if } Z=1
	\end{cases} $$ we have that its image by $B$ is	$B(h^*(Z)) = g(A)\left(\frac{h^*(0)}{h(0)}(1-P(Z=1|A)) + \frac{h^*(1)}{h(1)}P(Z=1|A) \right) $, which amounts to the following matrix times vector multiplication

	\[
	\begin{pmatrix}
		h_0^* \\
		h_1^*
	\end{pmatrix} \overset{B}{\longmapsto}
	\begin{pmatrix}
		\displaystyle\frac{(1 - g)(1 - b_0)}{1 - h} & \displaystyle\frac{(1 - g) b_0}{h} \\
		\displaystyle\frac{g(1 - b_1)}{1 - h} & \displaystyle\frac{g b_1}{h}
	\end{pmatrix}
	\begin{pmatrix}
		h_0^* \\
		h_1^*
	\end{pmatrix}
	\]
	
	with $P(Z^*=1):=h^*_1= 1-h^*_0$, $b_a:=P(Z=1|A=a)$, $g:=g(1)=P(A=1)$, $h=h(1)=P(Z=1)$ and $B$ represented by the matrix.

	 {For illustration, we set the natural policy to $h=0.6$ and $P(A=1\mid Z=z) = 0.3z + 0.5(1-z)$. Empirically, $\hat{b}_0 \simeq 0.68$ and $\hat{b}_1 \simeq 0.47$. For $g^*=(0.65,0.35)$, the unconstrained candidate $(B^*B)^{-1}B^*g^*$, i.e.\ the preimage of the projection of $g^*$ onto $B(\mathcal{F}_1)$, falls inside $\mathcal{F}_1$ and equals approximately $(0.25,0.75)$. However, for $g^*=(0.4,0.6)$ we obtain $(B^*B)^{-1}B^*g^* = (1.49,-0.49) \notin \mathcal{F}_1$, since this quantity is the minimizer of $\langle g^*-Bh^*, g^*-Bh^*\rangle$ over all of $\mathbb{R}^2$ rather than over $\mathcal{F}_1$. Projected gradient descent (Algorithm 2) is therefore required, and yields the constrained minimizer $(1,0)$, with corresponding implied intervention $g(h^*) = (0.5,0.5)$.}
	
\end{example}

Pick an arbitrary $g^*(A|W)\in \mathcal{G}_1 $. The candidate element $(B^*B)^{-1}B^*g^*$ to be the preimage of the least squares projection of $g^*$ onto $\mathcal{R}(B)$ is in general not in $\mathcal{F}_1$. This is happening because  $(B^*B)^{-1}B^*g^*$ is the solution to 
		$$\min_{h^* \in \mathcal{F}} \langle g^* -Bh^*,g^* - Bh^* \rangle_{L^2(P_{AW})}$$
	This is the same as finding $h^*$ such that $g^* - Bh^* \perp \operatorname{Ran}(B)$. Therefore we cannot look at projections on the rank. However, we pursue
	$$\min_{h^* \in \mathcal{F}_1} \langle g^* -Bh^*,g^* - Bh^* \rangle_{L^2(P_{AW})}$$
(mind the sets over which minimization takes place) so that our risk function is \begin{equation}\label{sq:risk-app}
	h^* \mapsto R_{g^*}(h^*) = \int (g^* - B(h))^2 dM(A,W)\end{equation}
	
	with $M=P_{AW}$. We arrived to a problem involving minimization of a convex functional over a convex set:

	\subsection{Least squares projection with continuous instrument}\label{shift}
	
Recall the definition of the indexing set $\mathcal{I}_n$ in \ref{sieve}. In contrast to containing densities, now the model will contain ratios so explicit presence of $h(Z|W)$ and the need to estimate it both vanish. We do not allow its dimension to grow with sample size and set it to $q$, which is the dimension of the instrument $Z$
	
$$\mathcal{I}_q= \left\{C(\beta)\exp{\left(\sum_{j} \beta_j \psi_j\right) : \beta(j) \in \mathbb{R}, j = 1,\ldots,q }\right\}$$
and we set its basis functions to be $$	\psi_j(Z,W): = -\frac{\partial}{\partial Z^j} \log h(Z \mid W)$$
Crucially, what belongs to $\mathcal{I}_q$ are not interventional densities anymore as in Section \ref{EM}, but ratios $h^*/h$. The following result clarifies that this particular choice of basis functions is motivated by being able to interpret the index $\beta$—which parameterizes an element $\frac{h^*(\beta)}{h} \in \mathcal{I}_q$—as representing the magnitude of a stochastic shift intervention \citep{diaz}.

\begin{propo}\label{clashift}
	Assume $h(\cdot |W)$ is differentiable $P(W)-$a.s. Let a stochastic intervention $h^*$ on the instrument be such that $\frac{h^*(\beta)}{h} \in \mathcal{I}_q$ with $	
	\psi_j(Z,W): = -\frac{\partial}{\partial Z^j} \log h(Z \mid W)$. Then $$h^*(\beta)(Z|W) = h(Z -\beta |W) + \mathcal{O}(\|\beta\|^2)$$
	\end{propo}
	
	We are just showing that apart from constituting an interesting class itself,  $\mathcal{I}_q$ locally approximates stochastic shift interventions with the particular choice of basis functions equal to the minus gradient of the natural density $h(Z|W)$ with respect to $Z$. Moreover, the $q$ dimensional	$\psi(Z,W): = \nabla_Z \log h(Z \mid W)$ is the score neural network of a diffusion model \citep{feng2024optimalconvexmestimationscore} so that it can be easily estimable without computing the whole natural density $h(Z \mid W)$.  An illustration of how well the basis functions are approximated can be found in Figure \ref{fig:scores} together with a basic explanation of score matching through diffusion models, both in Appendix \ref{app:diff}.

\begin{example}
	If $h(Z \mid W)$ is Gaussian with mean $\mu(W)$ and variance $\sigma^2(W)$, then:
	\[
	\psi(Z,W) = - \frac{Z - \mu(W)}{\sigma^2(W)}.
	\]
	The exponential tilt becomes
	\[
	\exp\left(-\beta(W)^T\frac{Z - \mu(W)}{\sigma^2(W)}\right)
	\]
	
	\noindent while the actual density ratio is (notice the $O(\beta(W)^2)$ term)
	\[
	\frac{h^*(Z \mid W)}{h(Z \mid W)} = \exp\left(-\beta(W)^T \frac{Z - \mu(W)}{\sigma^2(W)} - \frac{\beta(W)^2}{2 \sigma^2(W)}\right),
	\]
	
\end{example}

Now the conditional expectation involved in our operator $B$ is
	\[\begin{aligned}
{E}\left[ \frac{h^*(Z \mid W)}{h(Z \mid W)} \bigg| A, W \right] &=  	{E}\left[C(W){\exp\left(\beta(W)^T \psi(Z,W)\right)}\bigg| A, W \right] \\ &= C(W){E}\left[{\exp\left(\beta(W)^T \psi(Z,W)\right)}\bigg| A, W \right]\end{aligned}
\]
	
 Once we can evaluate this conditional expectation, we can compute its multiplicative action over $g(\cdot|W)$ for a given $W$. We are interested in fitting the parameters $\beta$ in this family w.r.t. a desired $g^*$. Our theoretical risk becomes

\begin{equation}\label{th:risk}
	\beta \in \mathbb{R}^q \longmapsto {E}_{AW}\left(g^*(A|W) -g(A|W )C(W){E}\left[{\exp\left(\beta^T \psi(Z,W)\right)}\bigg| A, W \right]\right)^2
	\end{equation}
assuming $	\beta (W) = \beta$, where $q$ is the dimension of the instrument.

The minimization of the risk in \eqref{th:risk} is carried out in two stages. First, for a fixed value of $\beta$, we estimate the inner conditional expectation
by regressing the transformed outcome $\exp(\beta^\top \psi(Z,W))$ on $(A, W)$ using a HAL basis and \texttt{glmnet} with \texttt{family=‘cox’}. This step is repeated each time $\beta$ is updated. In the second stage, the outer risk in \eqref{th:risk}—which compares the target intervention $g^* $ to the implied one—is minimized with respect to $\beta$ using an optimization routine such as \texttt{optim}. This process is iterative: for each candidate $\beta$ proposed during the optimization, the corresponding conditional expectation is re-estimated, and the objective is recomputed until convergence to an optimal $\beta^*$, which will correspond with how much do we have to shift the instrument as $Z^* = Z + \beta^*$ in order to be $L^2(P_{AW})$-closest to a chosen desired intervention $g^*$, according to Proposition \ref{clashift}.

\begin{figure}[t!]
	\centering
	\includegraphics[width=\textwidth]{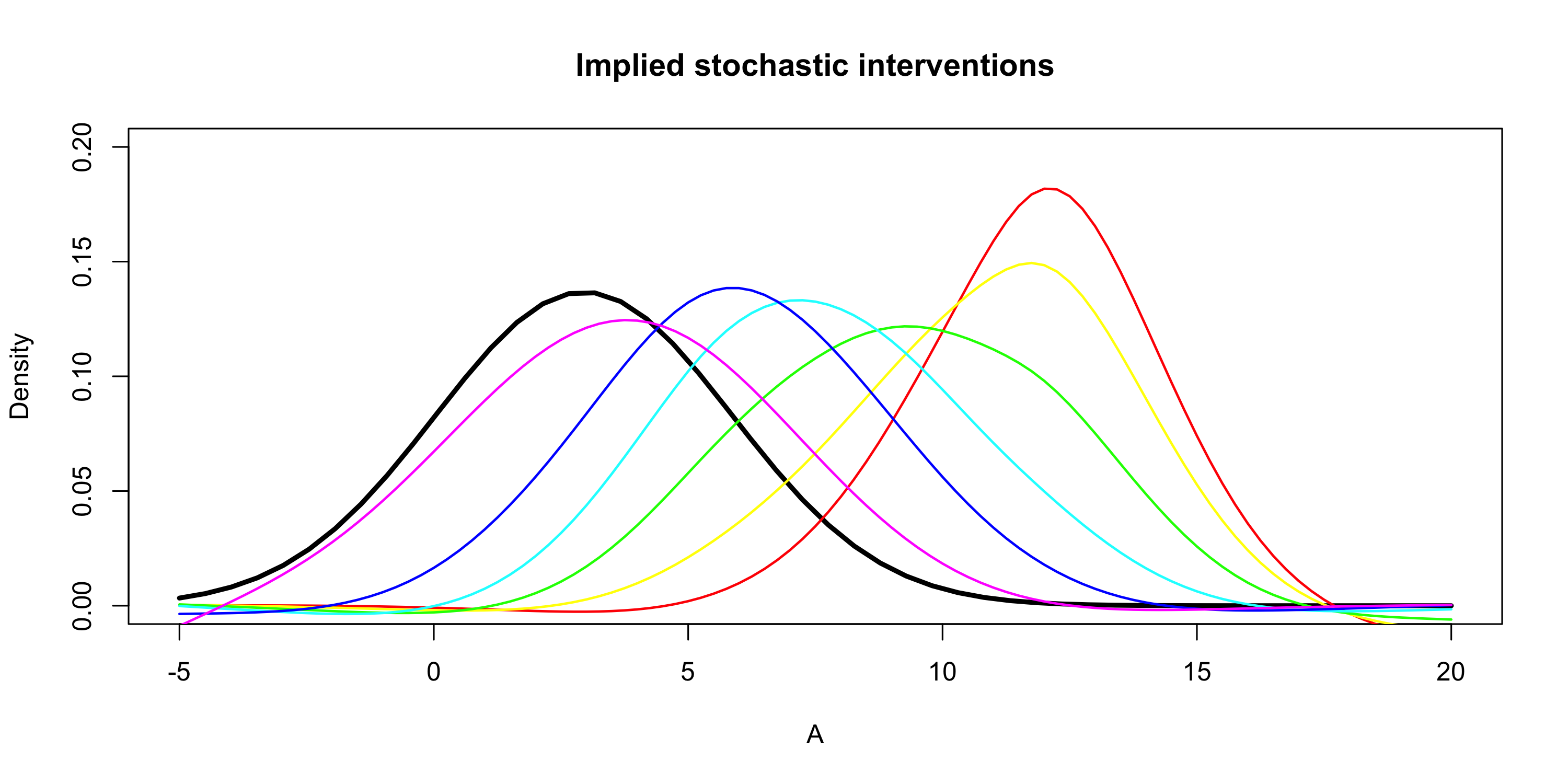}
	\caption{Implied interventional propensity scores along $\beta_2=-\beta_1$ and varying $\beta_1$ on a grid of $7$ values (colors) from $0.1$ (pink) to $2$ (red).  Recall $\beta$ represents the direction and magnitude by which the instrument is shifted—i.e., $Z^*  = Z + \beta$—in order for the resulting implied interventions to be the ones in the plot. The solid black line is the marginal observational distribution of $A$}
	\label{fig:implied}
\end{figure}

\begin{figure}[h!]
	\centering
	\includegraphics[width=\textwidth]{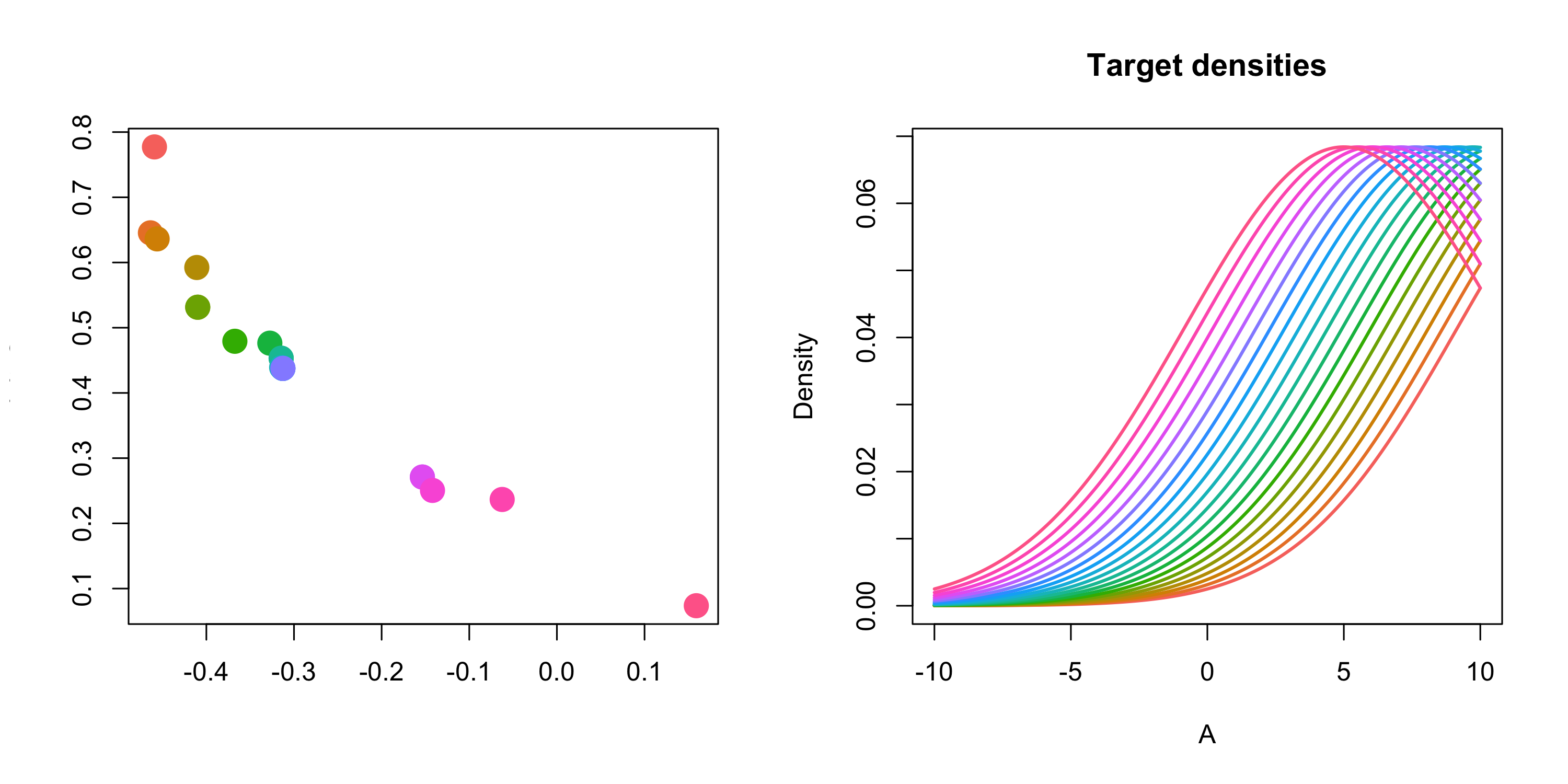}
	\caption{Each target density $g^*$ is colored according to its optimal $\beta^* = (\beta^{*1},\beta^{*2})$ corresponding to horizontal and vetical axes on the left plot obtained through empirical risk minimization of (\ref{th:risk}). Recall that the points on the left can be interpreted as bidimensional shifts of $Z$ in viritue of Proposition \ref{clashift}. }
	\label{fig:scatter}
\end{figure}

To illustrate Section \ref{shift}, we now present a simulation study with a bivariate continuous instrument $Z \in \mathbb{R}^2$ that demonstrates how score functions $\psi(Z,W)$ allow to express shift interventions on the instrument, how they propagate towards their corresponding implied interventions, and how effectively the minimization of \eqref{th:risk} recovers the best approximating implied intervention. We perform score matching—that is, estimation of $\psi_1(Z,W)$ and $\psi_2(Z,W)$—using a diffusion-based neural network implemented in \texttt{torch} with automatic differentiation (see Appendix \ref{app:diff}). We then examine two related tasks. First, we vary $\beta$ and compute the corresponding implied interventions indexed by $\beta$, visualized in Figure \ref{fig:implied}. Second, we fix 20 distinct target propensity scores $g^*$ and, for each, solve the projection problem in \eqref{th:risk} to find the optimal shift parameter $\beta^*$. These results are shown in Figure \ref{fig:scatter}.

The simulation is conducted with a sample size of $n=200$. Covariates $W = (W_1, W_2) \in \mathbb{R}^2$ are sampled independently from a uniform distribution: $W \sim \text{Uniform}([-2, 2]^2)$. We define a nonlinear conditional mean function $\mu(W) = (W_1,\ W_2^2)$. The instrument $Z \in \mathbb{R}^2$ is then generated conditionally on $W$ as $Z \mid W \sim \mathcal{N}(\mu(W), \sigma^2 I_2)$
with $\sigma$ a fixed noise scale. Finally, the treatment is generated as $A = W^\top \beta_W + Z^\top \beta_Z + \varepsilon $,
where $\beta_W = (0.3, 0.1)$, $\beta_Z = (3, -5)$, and $\varepsilon \sim \mathcal{N}(0, 0.1^2)$. This simulation design allows us to validate the accuracy of the estimated score function and to assess the efficacy of our projection method in approximating desired treatment policies via stochastic shifts.

\subsubsection{Sufficiency and the group of stochastic interventions}

We have tacitly reached the conclusion in Section \ref{shift} that stochastic interventions $h^*$ on the instrument whose ratio $h^*/h$ belongs to ${I}_q$ define a multiplicative group action 
	\[
	\begin{aligned}
		&G \times \mathcal{G}_1 \longrightarrow \mathcal{G}_1 \\
		& (\beta,g) \longmapsto g(A|W)C(W){E}\left[{\exp\left(\beta^T \psi(Z,W)\right)}\bigg| A, W \right]\
	\end{aligned}
	\]

	where $\mathcal{I}_q$ can be identified with a group $G$ of (local) shift interventions on the instrument $G$, meaning that the set of implied interventions $g(h^*)$ is the orbit of the natural $g$. 
	
	This structure of the intervened propensity score admits a natural analogy with exponential families. Recall that a probability density or mass function $p_\theta(Z)$ belongs to an exponential family if it can be written in the form
	\[
	p_\theta(Z) = h(Z) \exp\left( \eta(\theta)^\top \psi(Z) - A(\theta) \right),
	\]
	where $\psi(Z)$ and $h(Z)$ do not depend on $\theta$, $\eta(\theta)$ is the natural parameter and $A(\theta)$ is the log-partition function ensuring normalization. In this classical setting, the statistic $\psi(Z)$ is said to be sufficient for the parameter $\theta$ in the sense of Neyman:  $\psi$ is sufficient for $\theta$ if and only if one can factorize $p_\theta$ in the form $p_\theta(Z)=g_\theta(\psi(Z)) h(Z)$ for all $Z, \theta$ for some functions $g_\theta(\cdot) \geq 0$ and $h(\cdot) \geq 0$.
	
Our setting is a counterfactual analog of the Neyman factorization, where the intervened density factors into a baseline term and an intervention-dependent term involving $\psi$. Interventional ratios $ h^*/h \in \mathcal{I}_q$ are of the form:
	\[
	g(h^*)(A \mid W) = g(A \mid W) C(W)  {E}\left[ \exp\left( -\beta^\top \psi(Z, W) \right) \Big| A, W \right],
	\]
 This formulation resembles an exponential family in that it expresses the intervened density as a tilting of the baseline density through an exponential function of a statistic $\psi(Z, W)$. In this sense, we can say that the score $\psi(Z, W)$ involved in the backward denosing process of a diffusion model is acting as a \textit{causal }sufficient statistic.

\end{document}